\documentclass[12pt,thmsa]{amsart}
\def\Cal{\mathcal}

\def\F{{\Cal F}}

\def\T{{\Cal T}}
\def\Q{{\Cal Q}}

 \def\RR{{\mathbb R}}

\def\brm{\bbr^m}

\def\irm{\int_{\bbr^m}}

\def\irm1{\int_{\bbr^{m-1}}}

\def\bbr{{\Bbb R}}

\def\bbp{{\Bbb P}}
\def\bbn{{\Bbb N}}
\def\bbh{{\Bbb H}}

\def\bbc{{\Bbb C}}
\def\bbz{{\Bbb Z}}
\def\bbd{{\Bbb D}}

\def\sgn{{\hbox{\rm sgn}}}

\def\const{{\mbox{\rm const}}}
\def\cos{{\hbox{\rm cos}}}

\def\min{{\hbox{\rm min}}}

\def\dist{{\hbox{\rm dist}}}

\def\part{\partial}
\def\intl{\int\limits}
\def\b{\beta}
\def\vare{\varepsilon}

\def\Gam{\Gamma}
\def\Om{\Omega}
\def\a{\alpha}
\def\om{\omega}

\def\Del{\Delta}
\def\del{\delta}
\def\vp{\varphi}

\def\gam{\gamma}
\def\Lam{\Lambda}
\def\sig{\sigma}
\def\lam{\lambda}
\def\z{\zeta}
\def\e{\varepsilon}
\def\t{\tau}

\def\sideremark#1{\ifvmode\leavevmode\fi\vadjust{\vbox to0pt{\vss
 \hbox to 0pt{\hskip\hsize\hskip1em
\vbox{\hsize2cm\tiny\raggedright\pretolerance10000
 \noindent #1\hfill}\hss}\vbox to8pt{\vfil}\vss}}}%

\newtheorem{theorem}{Theorem}[section]
\newtheorem{lemma}[theorem]{Lemma}
\newtheorem{proposition}[theorem]{Proposition}

\theoremstyle{definition}
\newtheorem{definition}[theorem]{Definition}

\newtheorem{example}[theorem]{Example}

\theoremstyle{corollary}
\newtheorem{corollary}[theorem]{Corollary}
\theoremstyle{remark}
\newtheorem{remark}[theorem]{Remark}

\numberwithin{equation}{section}

\newcommand{\be}{\begin{equation}}
\newcommand{\ee}{\end{equation}}

\newcommand{\bea}{\begin{eqnarray}}
\newcommand{\eea}{\end{eqnarray}}
\newcommand{\Bea}{\begin{eqnarray*}}
\newcommand{\Eea}{\end{eqnarray*}}

\begin{document}

\title [The  Radon Transform ]{
 The  Radon Transform on the Heisenberg Group and the Transversal
Radon Transform}

\author{B. Rubin}
\address{Department of Mathematics, Louisiana State University, Baton Rouge,
Louisiana 70803, USA}
\email{borisr@math.lsu.edu}

\thanks{ The author was
supported  by the NSF grants PFUND-137 (Louisiana Board of Regents)
and DMS-0556157.}


\date{October 12, 2009}



\begin{abstract}

The notion of the Radon transform  $R_H$  on the Heisenberg group
$\bbh_n$ was introduced by R. Strichartz and inspired by  D. Geller
and E.M. Stein's  related work.
 The
 transversal Radon transform  $R_T$  integrates functions on $\bbr^m$,  $m\ge 2$,
 over hyperplanes meeting the last coordinate axis.
 In fact,  $R_H$ is a particular case of $R_T$, corresponding to $m=2n+1$.
We obtain new boundedness results and explicit inversion
 formulas for both transforms on $L^p$ functions  in the full range
 of the parameter $p$. We also show that  $R_H$ and $R_T$ are isomorphisms of the corresponding Semyanistyi-Lizorkin
spaces of smooth functions. In the framework of these spaces we obtain
inversion formulas, which are pointwise analogues of the corresponding formulas by R. Strichartz.

\end{abstract}

\maketitle

\section{Introduction}

\setcounter{equation}{0}

Let $\bbh_n=\bbc^n \times \bbr$ be the Heisenberg
 group with the multiplication law
 $$ (z,t) \circ (\z, \t)=(z+\z, t+\t-\frac{1}{2} \,Im (z\cdot \bar \z))$$
 and let
 \be\label{SRT} (R_H f)(z,t)=\int_{\bbc^n} f((z,t) \circ (\z, 0))\,
 d\z=\int_{\bbc^n} f(\z, t-\frac{1}{2} \,Im (z\cdot \bar \z))\,
 d\z\ee
be the  Radon-like transform, which was introduced by Strichartz in
his remarkable paper \cite{Str2}. He called $R_H$ the Heisenberg
Radon transform and noted that  definition (\ref{SRT}) stems from
the previous work by Geller and Stein \cite{GS1, GS2}. The
Heisenberg Radon transform represents an interesting object from the
point of view of both harmonic analysis and integral geometry.
Related transforms in the more general context of nilpotent Lie
groups and Siegel domains were studied by Felix \cite { Fe1, Fe2},
He and Liu \cite{HL1, HL2}, Peng and Zhang \cite { PZ}; see, also
\cite {He1, He2, LH, NT}.

A  variety of deep results in the area of harmonic analysis on the
Heisenberg
 group
 can be found in fundamental works by M.L. Agranovsky, C. A. Berenstein, D.C.
 Chang,  M. Cowling, G.B. Folland,
D. Geller,  A.
 Kor\'anyi, D. M\"uller,  D.H. Phong,  F. Ricci, E.M.
 Stein, S.  Thangavelu, J. Tie,
 and many other authors; see, e.g., \cite {BCT, CK,  KR, Kor, MRS, St2, Than}, and references therein.

Owing to  \cite[p. 386, 387]{Str2}, in a slightly different notation
we have the mixed norm estimate \be \label{mne0}
\left[\int_{\bbc^{n}}\left(\int_{\bbr}|(R_Hf)(z,t)|^r\,dt\right)^{p'/r}\,d\z\right]^{1/p'}\le
c\,\|f\|_p\,, \ee  $$1\le p\le p_0, \quad p_0\!=\!1+1/(2n+1),\quad
r=p/(2n+1-2np), \quad 1/p+ 1/p' =1.$$ The paper  \cite{Str2} (see
pp. 386, 399) also contains  the
 inversion formulas \be \label{Siv}
 R_H^{-1}=(4\pi)^{-2n} \Big
(\frac{\partial }{\partial t}\Big )^n\, R_H\,\Big (\frac{\partial
}{\partial t}\Big )^n,\qquad R_H^{-1}=\F^{-1} J \F_2, \ee where
$\F_2$ denotes the Fourier transform in the $t$-variable, $\F$ is
the Fourier transform in all variables, and
$(Jf)(x,y)=f(-2y^{-1}x,y)$.

 Our research is motivated by the following questions, that stem from [Str2] and remained open in  afore-mentioned publications.

 {\bf Question 1.} Is $R_H f$  well defined a.e. for $f\in L^{p} (\bbh_n)$,  when
 $p>p_0$? If yes, then what can we say about the boundedness of $R_H
 $ in this case?

 {\bf Question 2.}   What is the substitute of (\ref{Siv}) for $f\in L^{p} (\bbh_n)$?

 {\bf Question 3.}  To what class of functions are formulas (\ref{Siv})  applicable pointwise?

 In the present article we answer these questions in a more general context of the
 so-called
 transversal Radon transform $R_T$, which includes $R_H$ as a particular
 case.  Our method differs  in principle from those in the previous related works \cite{Str2, Fe1, Fe2, HL1, HL2,
 PZ}. To be more specific, let us rewrite (\ref{SRT}) by setting $z=u+iv, \z=\xi+i\eta$, where
 $u,v, \xi, \eta \in \bbr^n$. We obtain
$$(R_H f)(u +iv,t)=\int_{\bbr^{2n}} f(\xi, \eta,  t-\frac{1}{2}
(v\cdot \xi-u\cdot \eta))\, d\xi d\eta. $$
 Then we identify $f(\xi+i\eta, \t)$ with  a function $f(x)$ on $\bbr^m$,
$m\!=\!2n\!+\!1$, where \be\label{ccon1}  x\!\equiv \!(x',
x_m)\!\in\! \bbr^m,\qquad x'\!=\!(\xi, \eta)\!\in\! \bbr^{m-1},\quad
x_m=\t\in \bbr,\ee and set \be\label{ccon2} a=\frac{1}{2}(-v,u)\in
\bbr^{m-1}, \qquad b=t\in\bbr.\ee  This gives \be \label
{ccon1a} (R_H f)(u+iv,t)=(R_T f)(a,b),\ee where \be
\label {PRT} (R_T f)(a,b)=\int_{\bbr^{m-1}} f(x', a\cdot x'+b)\,dx'.
\ee Integral (\ref{PRT}) can be regarded as a  Radon transform
 of a
function $f$ on $\bbr^m$ associated with the hyperplane
\be\label{hp1} h=\{x=(x', x_m) \in \bbr^m: \, x_m=a\cdot x'+b\},\ee
 which is transversal to the last variable. Following Strichartz [Str2, p. 385], we call (\ref{PRT}) the {\it
transversal Radon transform}.
 It resembles a special case of the affine Radon transform
studied in  Gelfand's school \cite{GGG, GGV} and the parametric
Radon transform in \cite {Ehr}.  To distinguish notation, we write
\be (a,b) \in \tilde \bbr^m=\bbr^{m-1}\times \bbr.\nonumber\ee Note
that $m$ in (\ref{PRT}) may be an arbitrary integer $ \ge 2$. The
case of $m$ odd agrees with the Heisenberg Radon transform $R_H$.
Thus, basic results for $R_H $ will follow from those
for $R_T$.

\noindent{\bf Main results.} We   invoke the classical Radon
transform \be \label {RRT1} (Rf)(h)=\int_h f(x) \,d_h (x), \ee where
integration over the hyperplane $h$ in $\bbr^m$ is performed against
the usual Lebesgue measure. Suppose that   $h$ is parameterized as
\be\label{hp2} h=\{x \in \bbr^m: \, x\cdot \theta=t\}, \qquad
(\theta,t) \in \bbp^m=S^{m-1}\times \bbr,\ee
   $S^{m-1}$ being the unit sphere in $\bbr^m$. Then
   \be \label {RRT} (Rf)(h)\equiv (Rf)(\theta,t)=\int_{\theta^\perp} f(y+t\theta)\,dy. \ee
    The corresponding
   dual transform is  \be\label{durt}
(R^*\psi)(x)= \int_{S^{m-1}} \psi(\theta, x\cdot \theta)\,d\theta.
\ee  The  duality relation  \be\label{duar}
\int_{\bbp^m} (Rf)(\theta, t) \,\psi(\theta, t)\,d\theta\, dt=
\int_{\bbr^m} f(x)(R^* \psi)(x) \,dx\nonumber\ee holds provided that at
least one of these integrals is finite when $f$ and $\psi$ are
replaced by $|f|$ and $|\psi|$, respectively.

Integral transforms (\ref{RRT1}) and (\ref{durt}),
   are the most familiar to the reader. They
were introduced by Radon \cite{Rad} and studied in numerous
publications together with their $k$-plane generalizations; see,
e.g., \cite{Ber, Dea, GGG, Hel, Kat, Mar, Na, Pal, QCK, RK, Ru7,
SSW} and references therein.

Unlike (\ref{hp2}), that parameterizes {\it all}  hyperplanes in
$\bbr^m$,  (\ref{hp1}) excludes  hyperplanes parallel to the
$x_m$-axis.  If both parameterizations are available, then the
standard  Calculus yields \be\label{con1} (Rf)(h)=\sqrt{1+|a|^2}\,
(R_T f)(a,b).\ee

Equation (\ref{con1}) means that basic properties of $R_H$ and $R_T$
can be derived from known facts for $R$. However, as indicated by
Strichartz \cite [p. 385]{Str2},  there are  essential distinctions
between $R_T$ and $R$. This remark inspires us to give direct
proofs, that might be instructive.

\begin{theorem}\label{t1} Let $Re  \,\a >0$,
\be
 \lam\!=\!\frac{\pi^{(m-1)/2} \, \Gam (\a/2)}{\Gam
((\a\!+\!m\!-\!1)/2)}.\nonumber\ee The following equalities hold:
\be \label{eq1}\int_{\tilde \bbr^m} (R_T f)(a,b) \frac{|b|^{\a-1}\,
da\, db}{(1+|a|^2)^{(\a+m-1)/2}} =\lam \, \int_{\brm}
f(x)\,|x|^{\a-1} \,dx, \ee \be \label{eq2}\int_{\tilde \bbr^m} (R_T
f)(a,b) \frac{|b|^{\a-1}\, da\, db}{(1+|a|^2+b^2)^{(\a+m-1)/2}}
=\lam \int_{\brm} f(x)\frac{|x|^{\a-1}}{(1+|x|^2)^{\a/2}} \,dx, \ee
provided that either side of the corresponding equality is finite
when $f$ is replaced by $|f|$.
\end{theorem}

Equalities (\ref {eq1}) and (\ref {eq2}) are proved in Section 2 and
play a key role. As we shall see in Sections 3 and 4,
 (\ref {eq1}) paves the way to the ``right" definition of the {\it
dual transversal Radon transform}, the corresponding analytic
families of fractional integrals of the Semyanistyi type, and a
variety of explicit inversion formulas for $R_T$. Equality (\ref
{eq2}) answers the question about existence of the transversal Radon
transform of $L^p$ functions and provides some boundedness results.
In particular, by H\"older's inequality,
 (\ref {eq2}) implies the following corollary (cf. \cite [Corollary 2.4]{Ru7}) .
\begin{corollary} \label{cco} If $f \in L^p(\brm), \; 1 \le p < m/(m-1)$, then $(R_T f)(a,b)$ is finite
for almost all $(a,b) \in \tilde \bbr^m$. Moreover, for any
$\a>1-m/p'$, $1/p+1/p'=1$, we have \be\label{est1} \int_{\tilde
\bbr^m} |(R_T f)(a,b)| \frac{|b|^{\a-1}\, da\,
db}{(1+|a|^2+b^2)^{(\a+m-1)/2}}\le c\, ||f||_{L^p
(\brm)},\nonumber\ee
 where $c=c(\a,p,m)= const < \infty$.
 \end{corollary}

Corollary \ref {cco} leads to the following result for the
Heisenberg Radon transform (just set $m=2n+1$).
\begin{theorem}\label{t2}  Let $p_1=1+1/2n$.
If $f \in L^p(\bbh_n), \; 1 \le p < p_1$, then $(R_H f)(z,t)$ is
finite for almost all $(z,t) \in \bbh_n$. Moreover, for any
$\a>1-(2n+1)/p'$, \be\label{est2} \int_{\bbh_n} |(R_H
f)(z,t)| \frac{|t|^{\a-1}\, dz\, dt}{(1+|z|^2+t^2)^{(\a+2n)/2}}\le
c\, ||f||_{L^p (\bbh_n)}.\nonumber\ee
 \end{theorem}

Clearly,  $p_1\!=\!1\!+\!1/2n$ in Theorem \ref {t2} is greater than
$p_0\!=\!1\!+\!1/(2n\!+\!1)$ in  (\ref{mne0}). This gives the
answer to Question 1. Note also that the condition $p < m/(m-1)$ in
Corollary \ref{cco} (resp. $p < p_1$ in Theorem \ref{t2}) is sharp.
This follows from (\ref{con1}) and the well-known fact that if $p
\ge m/(m-1)$ and $$ f(x)=(2+|x|)^{-m/p} (\log (2+|x|))^{-1} \; (\in
L^p),$$
 then $(Rf)(h) \equiv \infty$; see, e.g., \cite{So, Ru7}.

\begin{theorem}\label{tos} For $m\ge 2$ an
a priori
 inequality
$$
\left(\, \int_{\bbr^{m-1}} \left[\int_\bbr | (R_Tf) (a,b)|^r
db\right]^{q/r} da\right)^{1/q} \le c_{p, q, r}  \| f \|_p$$ holds
if and only if
\be\label{TRTaz}1\le p < m/(m-1), \qquad q
=p^{\prime},\qquad 1/r = 1-m/p'.\ee
\end{theorem}

This statement will be proved in Section 8. For $m=2n+1$, it
extends Strichartz's result (\ref{mne0}) to all $1\le p\!< \!p_1$,
where the bound $p_1\!=\!1\!+\!1/2n\!> \! p_0$ is best possible.

A variety of inversion formulas for the transversal Radon transform
and the Heisenberg  Radon transform are presented by Theorems \ref
{441a}, \ref{invek2}, \ref{441ab}, \ref{cbp}, \ref {cbpx},
\ref{inveH1}, and
 \ref{inveH2}. We
develop severals methods involving hypersingular integrals, powers
of minus-Laplacian, and wavelet transforms. We also present the
convolution backprojection inversion method for $R_Tf$, $f\in L^p
(\bbr^m)$, introduce the relevant ridgelet-like transforms and give
examples of admissible wavelet functions.  All inversion formulas for $L^p$-functions
are obtained in the maximal range of the parameter $p$.

Another series of pointwise inversion formulas  for $R_T$ and $R_H$
 is obtained in the framework of the spaces of the
Semyanistyi-Lizorkin  type.  For example, we introduce the space
$\Phi (\bbh_n)$,  which consists of Schwartz functions $\vp$ on
$\bbh_n$ satisfying
 \be  \int_{\bbr} \phi(\z, \t)\, \t^{k}
 \,d\t\!=\!0 \qquad \forall k\in \bbz_+, \quad \forall \z \in\bbc^{n}.\nonumber\ee

We show that the Heisenberg  Radon transform $R_H$ is  an
automorphism of $\Phi (\bbh_n)$ and obtain pointwise analogues of
 formulas (\ref{Siv}). In particular, by Theorem \ref{inve2} a
function $f \in \Phi (\bbh_n)$ can be reconstructed from $\vp=R_H f$
by the formula \be \label {ansiv}f
(\z,\t)=(-1)^n\,(4\pi)^{-2n}\,(\partial_\t^{n}R_H
\partial_t^{n}\vp)(\z,\t),\ee
which differs from the first  equality in
(\ref{Siv}) by the factor $(-1)^n$. To obtain these results we first establish more general
 relations for the  transversal Radon transforms, including Riesz potentials.

Some historical comments are in order. In 1960, Semyanistyi
\cite{Se1} came up with an interesting idea to
 treat the hyperplane Radon transform and its dual
 as members of suitable analytic families of operators $R^\a$ and $\stackrel{*}{R}{}\!\!^\a$,
 $\a \in \bbc$,  so  that  $R^\a f|_{\a =0}=Rf$
and $\stackrel{*}{R}\!{}^{\!\!\a}\vp|_{\a =0}=R^*\vp$. Semyanistyi's
operators combine properties of Radon transforms and fractional
integrals. They yield a variety of inversion formulas
$$ c_{n}f=(-\Delta)^{(n+\a-1)/2} \stackrel{*}{R}{}\!^\a (Rf),$$
where $c_n=\const$ and $(-\Delta)^{(n+\a-1)/2}$ is a power of
minus-Laplacian, that can be realized in several ways. This
philosophy was extended in different directions and proved to be
useful in
 applications to
convex geometry; see \cite {OR, Ru2}, \cite {Ru4}-\cite {Ru8},
\cite{ RZ}.

In the present article  we extend Semyanistyi's idea to the
transversal Radon transform and demonstrate its benefits.

We expect that  some ideas  of
 our article are applicable to  Radon transforms on more general two-step
 nilpotent groups, as in \cite {Fe1, Fe2, PZ, HL1, HL2}.
We plan to address this topic in forthcoming publications.

\noindent  {\bf Plan of the paper.}   Section 2 contains
 preliminaries  and the proof of Theorem \ref{t1}.  In Section 3 we
 introduce analytic families  of fractional integrals of the
 Semyanistyi type,
 associated to the transversal Radon transform $R_T$, and define the dual transversal Radon
 transform. This section forms a foundation for different
 inversion methods developed in Section 4.
 In Section 5 we introduce and investigate a function  space $\Phi (\bbr^m)$ of the
Semyanistyi-Lizorkin  type associated to the transversal Radon
transform. Section 6 is devoted to the detailed study of the
 connection between the classical hyperplane Radon
transform (\ref{RRT1}), the
 transversal Radon transform (\ref{PRT}), the Heisenberg Radon
transform (\ref{SRT}), and their duals. Results of this section are
applied in Sections 7 to investigation of isomorphism properties and
inversion of the Heisenberg  Radon transform.
 Section 8 deals with mixed norm estimates of the
transversal and Heisenberg  Radon transforms.

\noindent   {\bf Acknowledgements}. The work was started when I
 was visiting the  Capital Normal University
 (Beijing) and the Guangzhou University   in Summer 2009. I am  grateful to  Professors  Zhongkai Li and Jianxun He
  for the hospitality, support, and discussions. They participated
 in some calculations related to  Sections 5-8, and their results will be published elsewhere. I am especially  thankful
 to Professor Dr. Rainer Felix for sending me his interesting paper \cite
 {Fe2}.

\section{ Preliminaries. Proof of Theorem \ref{t1}}

\setcounter{equation}{0}

\subsection{Notation and auxiliary statements}
We denote by  $e_1, \ldots, e_m$ the coordinate unit vectors in
$\brm$; $\theta=(\theta_1, \ldots, \theta_m)\in S^{m-1}$;
$\;\theta'=(\theta_1, \ldots \theta_{m-1})$;
 $$S^{m-1}_+=\{\theta\in
S^{m-1}:\,\theta_m>0\};$$ $\sig_{m-1}=2\pi^{m/2}/\Gam (m/2)$ is the
area of $S^{m-1}$; $d\theta$ stands for the usual Lebesgue measure
on $S^{m-1}$. We recall the well-known
 Catalan formula from Calculus
\be \label{hva37} \int_{S^{m-1}} f( \theta \cdot \sig)\,
d\theta=\sig_{m-2} \int_{-1}^1 f(t) (1-t^2)^{(m-3)/2}\, dt,\ee which
holds for every $ \sig\in S^{m-1}$ provided  that at least one of
the integrals in (\ref{hva37}) is finite when $f$ is replaced by
$|f|$; see, e.g., \cite[p. 216]{ABR}. The notation $C$,  $C^\infty$,
and $L^p$ for function spaces  on
 $\bbr^m$ is classical; $C_0=\{f\in C(\bbr^m):
\lim\limits_{|x|\to\infty} f(x) = 0\}$; $S(\bbr^m)$ is the Schwartz
space of rapidly decreasing $C^\infty$-functions  with standard
topology. All function spaces on the Heisenberg
 group $\bbh_n=\bbc^n \times \bbr$ are identified with the
 corresponding spaces on $\bbr^{2n+1}$.

The Riesz potential $I^\a f$ on $\bbr^m$ is defined by
\be \label{15}(I^\a f)(x)=\frac{1}{\gamma_m(\a)} \int_{\bbr^m}
 \frac{f(y)\,dy}{|x-y|^{m-\a}},\quad
  \gamma_{m}(\a)=
  \frac{2^\a\pi^{m/2}\Gamma(\a/2)}{\Gamma((m-\a)/2)},
  \ee
$$Re \, \a >0, \qquad \a-m \neq 0,2,4, \ldots ;$$
see \cite{St1}. The letter $c$ stands for a constant that can be different at each
occurrence; $\bbz_+ =\{0,1,2, \ldots \}$; $\bbz_+^m$ denotes the
relevant
 set of  multi-indices. Given a real-valued expression
$A$, let $(A)_+^\lam=A^\lam$ if $A>0$ and $(A)_+^\lam=0$ if $A \le
0$.

\begin{lemma}\label{tlem}${}$\hfil

\noindent ${\rm {(i)}}$ If $f\in L^1(S^{m-1})$, then
 \be\label{teq1} \int_{S^{m-1}}
f (\theta)\, d\theta=\int_{\bbr^{m-1}}\tilde f (a)\,\tilde d
 a,\ee
$$
\tilde f (a)=f \Big (\frac{a+e_m}{\sqrt{1+|a|^2}} \Big )+f \Big
(\frac{a-e_m}{\sqrt{1+|a|^2}} \Big ),\qquad \tilde d
 a=\frac{da}{(1+|a|^2)^{m/2}}.$$

\noindent ${\rm {(ii)}}$ Conversely, if $f\in L^1(\bbr^{m-1})$, then
 \be\label{teq2} \int_{\bbr^{m-1}} f(a)\, da=\int_{S^{m-1}_+}f
\Big (\frac{\theta'}{\theta_m}\Big )\,
\frac{d\theta}{\theta_m^m}.\ee
\end{lemma}
\begin{proof}  We write
$$
\int_{S^{m-1}} f (\theta)\, d\theta=\int_0^\pi  \sin ^{m-2} \om\,
d\om\int_{S^{m-2}} f (\sig\sin\om +e_m \,\cos\,\om)\, d\sig$$ and
then set $s=\tan \om$. This gives
$$
\int_0^\infty \frac{s^{m-2}}{(1+s^2)^{m/2}}\int_{S^{m-2}}\left [ f
\Big (\frac{s\sig +e_m}{\sqrt{1+s^2}} \Big )+f \Big (\frac{s\sig
-e_m}{\sqrt{1+s^2}} \Big )\right ]\, d\sig,$$ which coincides with
 (\ref{teq1}). Similarly,
\bea
 &&\int_{\bbr^{m-1}} f(a)\, da=\int_0^\infty s^{m-2}\, ds\int_{S^{m-2}}
f (s\sig)\, d\sig\nonumber \\
&=&\int_0^{\pi/2}\frac{\sin ^{m-2} \om}{\cos^m \om }\,
d\om\int_{S^{m-2}}
 f \Big (\frac{\sig\sin\om}{\cos\, \om}\Big )\, d\sig=\int_{S^{m-1}_+}f
\Big (\frac{\theta'}{\theta_m}\Big )\,
\frac{d\theta}{\theta_m^m}.\nonumber\eea
 \end{proof}

We will be dealing with Riemann-Liouville fractional integrals
\be\label{ril} (I^\a_{0+} \psi)(\tau) =
\frac{1}{\Gamma(\a)}\int^\tau_0 (\tau-s)^{\a-1}\psi(s) \,ds.\ee The
following  fact from  \cite {Ru3} will be needed, where $[\a]$
denotes the integer part of $\a$.
\begin{lemma} \label {frai}  Let $ \a>0$,
\be\label{216} \int_0^\infty s^j\psi(s)\,ds=0 \ \mbox{ for all }
j=0,1,\ldots,[\a], \nonumber\ee \be\label{217} \int_0^\infty
s^\b|\psi(s)|\,ds<\infty \ \mbox{ for some } \b>\a.\nonumber\ee Then
\be\label{218} |(I_{0+}^{1+\a}\psi)(s)|=\left\{ \!
 \begin{array} {ll} \! O(s^{\a})&\mbox{ if $ 0<s\le 1$},\\
\!O(s^\gam),\;\gam=\a-\min(1+[\a],\b)<0,&\mbox{ if $s>1$}.\\
 \end{array}
\right. \nonumber\ee Furthermore, \be\label{219}
\int_0^\infty(I_{0+}^{1+\a}\psi)(s)\frac{ds}{s}= \left\{ \!\!
 \begin{array} {ll}  \displaystyle{\Gam(-\a) \int_0^\infty s^\a\psi(s)\,ds} &\mbox{ if $\a\notin\bbn$},\\
{}\\ \!\displaystyle{\frac{(-1)^{\a+1}}{\a!}\int_0^\infty s^\a\psi(s)\log s\, ds} &\mbox{ if $ \a\in\bbn$}.\\
 \end{array}
\right.\nonumber\ee
\end{lemma}

  More information about Riesz potentials and fractional
integrals can be found in [Ru1], [SKM].

\subsection{Proof of Theorem \ref{t1}}  We denote by $I$ the left-hand
side of (\ref{eq1}). Then \bea
I&=&\int_{\bbr^{m-1}}\frac{da}{(1+|a|^2)^{(\a+m-1)/2}}\int_{\bbr}
|b|^{\a-1}\, db\int_{\bbr^{m-1}} f(x', a\cdot x' +b)\,
dx'\nonumber\\
&=&\int_{\bbr^{m}} f(x) \, dx\int_{\bbr^{m-1}}\frac{|a\cdot x'
 -x_m|^{\a-1}}{(1+|a|^2)^{(\a+m-1)/2}}\, da.\nonumber\eea
By (\ref{teq2}) and (\ref {hva37}), the inner integral can be represented as
$$
\int_{S^{m-1}_+} |\theta \cdot x|^{\a-1}\,
d\theta=\frac{|x|^{\a-1}}{2}\int_{S^{m-1}} |\theta \cdot
e_m|^{\a-1}\, d\theta=\frac{\pi^{(m-1)/2} \, \Gam (\a/2)}{\Gam
((\a+m-1)/2)}\, |x|^{\a-1},
$$
and (\ref{eq1}) follows.

The proof of (\ref{eq2}) uses the same idea. Denoting by $I$ the
left-hand side,  we get $$I=\int_{\bbr^{m}} f(x) \,
dx\int_{\bbr^{m-1}}\frac{|a\cdot x'
 -x_m|^{\a-1}\, da}{(1+|a|^2+|a\cdot x'
 -x_m|^2)^{(\a+m-1)/2}}. $$
 By (\ref{teq2}) and (\ref{hva37}), the inner integral becomes
$$
\frac{1}{2}\int_{S^{m-1}} \frac{|\theta \cdot x|^{\a-1}\,
d\theta}{(1+|\theta \cdot x|^2)^{(\a+m-1)/2}}=\sig_{m-2}
r^{\a-1}\int_0^1 \frac{t^{\a-1}
(1-t^2)^{(m-3)/2}}{(1+r^2t^2)^{(\a+m-1)/2}}\, dt,
$$
where $r=|x|$.  The last integral can be computed by changing
variables and using \cite[formula 3.238(3)]{GRy}. \hfill $\square$

\section{Analytic families of fractional integrals and the dual transversal Radon transform}
Below we introduce analytic families of fractional integrals
associated to the transversal Radon transform. We start with
(\ref{eq1}) and   replace $f$ by the shifted function
$f_x(y)=f(x+y)$, so that \be\label{sh} (R_T
f_x)(a,b)=(R_Tf)(a,b-a\cdot x' +x_m). \ee Changing variables,
 we get
\bea \label {32} &{}&\frac{1}{\Gam (\a/2)}\int_{\tilde \bbr^m} (R_T
f)(a,b) \frac{|a\cdot x'+b -x_m |^{\a-1}\, da\,
db}{(1+|a|^2)^{(\a+m-1)/2}}
\\
&=&\frac{\pi^{(m-1)/2}}{\Gam ((\a+m-1)/2)}
 \int_{\bbr^{m}}
  f(y)
|x-y|^{\a-1} \, dy, \quad Re \, \a >0. \nonumber \eea Using the
normalizing coefficient of the one-dimensional  Riesz potential (cf.
(\ref{15})), we denote \be\label{dt1} (\stackrel{*}{R_T}{}\!\!^\a
\vp)(x)=\frac{1}{\gamma_{1}(\a)} \int_{\tilde \bbr^m} \vp(a,b)\,
\left (\frac{|a\cdot x'+b -x_m |}{\sqrt{1+|a|^2}}\right )^{\a-1} \,
\tilde d a\, db,\ee
$$ x\in \bbr^m, \qquad Re \, \a >0, \qquad \a\neq 1,3,5, \ldots; \qquad \tilde d a=\frac{da}{(1+|a|^2)^{m/2}}.$$
Under the same assumption for $\a$, we  define \be
\label{dt2}(R_T^\a f)(a,b)=\frac{1}{\gamma_{1}(\a)}\int_{\bbr^m}
f(x) \left (\frac{|a\cdot x'+b -x_m |}{\sqrt{1+|a|^2}}\right
)^{\a-1} \, dx, \ee $(a,b) \in \tilde \bbr^m$. The following duality
relation is an immediate consequence of Fubini's theorem:
 \be \label{dual1} \int_{\tilde \bbr^m} (R_T^\a f)(a,b)\,
\vp(a,b)\,\tilde d a\, db=\int_{\bbr^m} f(x)\,
(\stackrel{*}{R_T}{}\!\!^\a \vp)(x)\, dx.\nonumber\ee Note that \be
\frac{|a\cdot x'+b -x_m |}{\sqrt{1+|a|^2}}=\dist (x,h) \nonumber\ee
is the Euclidean distance between the point $x\in \bbr^m$ and the
hyperplane $h=\{y\in \bbr^m: a\cdot y' +b=y_m\}$.
 Thus, the kernel of fractional integrals (\ref{dt1}) and (\ref{dt2})
 has a natural geometric meaning. By  Fubini's theorem,
 \be\label{rep1}
R_T^\a f=\psi_\a I_1^\a R_T f,\qquad  \stackrel{*}{R_T}{}\!\!^\a \vp
=\stackrel{*}{R_T} \psi_\a I_1^\a\vp,\ee where \be\psi_\a
(a)=(1+|a|^2)^{(1-\a)/2},\nonumber\ee $I_1^\a$ denotes the
one-dimensional Riesz potential of order $\a$ in the $b$-variable,
and \be\label{du1} (\stackrel{*}{R_T} \vp)(x)
=\int_{\bbr^{m-1}}\vp(a, x_m- a\cdot x')\, \tilde d  a. \ee We call
$(\stackrel{*}{R_T} \vp)(x)$ the {\it dual transversal Radon
transform} of $\vp$. By Fubini's theorem, \be \label{dual1}
\int_{\tilde \bbr^m} (R_T f)(a,b)\, \vp(a,b)\,\tilde d a\,
db=\int_{\bbr^m} f(x)\, (\stackrel{*}{R_T} \vp)(x)\, dx\ee provided
that at least one of these integrals is finite when $f$ and $\vp$
are replaced by $|f|$ and $|\vp|$, respectively. Indeed, \bea
l.h.s&=&\int_{\bbr^{m-1}}\tilde d a\int_{\bbr^{m-1}}dx'\int_\bbr
\vp(a,b)\,f(x', a\cdot x' +b)\, db\nonumber\\ &=&\int_{\bbr^m}
f(x)\, dx\int_{\bbr^{m-1}}\vp(a,x_m- a\cdot x')\,\tilde d
a=r.h.s.\nonumber\eea

 For sufficiently
good $f$ and $\vp$, owing to the limit property of  Riesz
potentials, we can pass to the limit in (\ref{rep1}) as $\a\to 0$.
This gives \be\label{311} \lim\limits_{\a \to 0 }R_T^\a f= \tilde
R_T f,\qquad \lim\limits_{\a \to 0
}\stackrel{*}{R_T}{}\!\!^\a\vp=\stackrel{*}{R_T} \tilde \vp,
\nonumber\ee
 $$(\tilde R_T f)(a,b)\!=\!\sqrt{1\!+\!|a|^2}\, (R_T f)(a,b),
\qquad \tilde \vp (a,b)\!=\!\sqrt{1\!+\!|a|^2}\,\vp (a,b). $$

  \begin{theorem}\label{t31} Let $\vp=R_Tf, \; f \in L^p(\bbr^m)$, \; $$1 \le p
  <\frac{m}{\a+m-1}, \qquad Re \, \a >0, \quad \a\neq 1,3,5, \ldots.$$
Then for almost all $x\in \bbr^m$, \be\label{33}
 (\stackrel{*}{R_T}{}\!\!^\a
\vp)(x)=(2\pi)^{m-1}\, (I^{\a+m-1} f)(x).\ee Moreover, for $1 \le p
  <m/(m-1)$ and almost all $x\in \bbr^m$,
  \be\label{33a}
 (\stackrel{*}{R_T}
\tilde \vp)(x)\!=\!(2\pi)^{m-1}\, (I^{m-1} f)(x),\quad \tilde \vp
(a,b)\!=\!\sqrt{1\!+\!|a|^2}\,\vp (a,b).\ee
\end{theorem}
\begin{proof} Equality (\ref{33}) coincides with (\ref{32}) above.
Equality (\ref{33a}) can be formally derived from (\ref{33}) by
passing to the limit as   $\a \to 0$. A rigorous proof for arbitrary
 $ f \in L^p(\bbr^m)$, \; $1 \le p
  <m/(m-1)$, is the following. Note that for such $p$ the Riesz potential
  $(I^{m-1} f)(x)$ is finite for almost all $x\in \bbr^m$.

  Owing to (\ref{sh}), it suffices to prove (\ref{33a}) at $x=0$. Thus,  we have to show that
  \be
  I\equiv \int_{\bbr^{m-1}} \frac{(R_T f)(a,0) \,da}{(1\!+\!|a|^2)^{(m-1)/2}}=
  \frac{\pi^{(m-1)/2}}{\Gam ((m-1)/2)} \, \int_{\bbr^{m}}\frac{f(y)}{|y|}\, dy.\nonumber\ee
  Changing the order of integration and passing to polar coordinates, we obtain
  $$
  I=\int_{S^{m-2}}d\sig \int_{\bbr^{m-1}} dy'\int_0^\infty
  \frac{r^{m-2}\, f(y', r\sig\cdot y')}{(1+r^2)^{(m-1)/2}}\, \, dr.$$
  Then, changing variables, we get
  $$
  I=\int_{\bbr^{m}} f(y) |y_m|^{m-2}\,J(y)\, dy,
  \qquad J(y)=\frac{1}{2}\,\int_{S^{m-2}}\frac{d\sig }{(y_m^2+(\sig\cdot y')^2)^{(m-1)/2}}.$$
  If $m=2$, then   $$J(y)=\frac{1}{2}\,\left [\frac{1}{(y_2^2+(1\cdot y_1)^2)^{1/2}} +
  \frac{1}{(y_2^2+(-1\cdot y_1)^2)^{1/2}}\right ] =\frac{1}{|y|}.$$
  If $m>2$, then $J(y)$ can be evaluated using   (\ref {hva37})
  and \cite[formula 3.238(3)]{GRy} as follows. Set $\lam=|y_m|/|y'|, \; \om =y'/|y'|$. Then
  \bea
  J(y)&=&\frac{1}{2|y'|^{m-1}}\int_{S^{m-2}}\frac{d\sig }{(\lam^2+(\sig\cdot \om)^2)^{(m-1)/2}}\nonumber\\&=&
\frac{\sig_{m-3}}{|y'|^{m-1}}\int_0^1\frac{(1-t^2)^{(m-4)/2}\,
dt}{(\lam^2+ t^2)^{(m-1)/2}}=\frac{\pi^{(m-1)/2}}{\Gam((m-1)/2)}\,
\frac{|y_m|^{2-m}}{|y|}.\nonumber\eea This gives what we need.
\end{proof}

\section{Inversion of the Transversal Radon Transform}

By Theorem \ref {t31}, inversion of the transversal Radon transform
$R_T$ reduces to inversion of Riesz potentials of $L^p$-functions on
$\bbr^m$. The last topic was studied in numerous publications; see
\cite{SKM, Ru1, Sam3} and references therein. The
 corresponding results are adapted in \cite {Ru7} for the $L^p$-theory
 of the $k$-plane Radon transforms on $\bbr^n$.

Theorem \ref {t31} and the second equality in (\ref{rep1}) suggest
 the following  approaches to reconstruction of $f$ from $\vp=R_T f$.

\vskip 0.3 truecm

\noindent {\bf The first approach:} \hfill

\noindent Invert the $m$-dimensional Riesz potential $I^{\a+m-1} f$
in (\ref {33}) to get \be\label{opt1}
 f=(2\pi)^{1-m} \, \bbd^{\a+m-1} \stackrel{*}{R_T}{}\!\!^\a\vp,\ee
where $\bbd^{\a+m-1}=(-\Del)^{(\a+m-1)/2}$ is the Riesz fractional
derivative. Here $\a$ can be chosen as we please.

\vskip 0.3 truecm

\noindent{\bf The second approach:} \hfill

\noindent Set (formally) $\a=1-m$ in (\ref {33}) to get
\be\label{opt2}f=(2\pi)^{1-m} \,
\stackrel{*}{R_T}{}\!\!^{1-m}\vp.\ee

Of course,
 $\stackrel{*}{R_T}{}\!\!^{1-m}\vp$ and all Riesz
derivatives should be properly interpreted.

The described  method is known in the literature as the method of
Riesz potentials; cf. \cite [Section 4]{Ru7}, \cite{Rou}. Below we
 provide formulas (\ref{opt1}) and (\ref{opt2}) with precise meaning.


\subsection{The first approach}

\subsubsection{ Hypersingular integrals} We assume $\a=0$, when (\ref
{33a}) yields  \be\label{opt1a} f=(2\pi)^{1-m} \,
\bbd^{m-1}\stackrel{*}{R_T} \tilde \vp, \qquad \tilde \vp
(a,b)\!=\!\sqrt{1\!+\!|a|^2}\,(R_T f) (a,b).\nonumber\ee Consider
finite differences \bea \label{diff1} &&(\Del^\ell_y
g)(x)=\sum_{j=0}^\ell
{\ell \choose j} (-1)^j g(x-jy), \\
\label{diff2} &&(\tilde\Del^k_y g)(x)=\sum_{j=0}^k {k \choose j}
(-1)^j g(x-\sqrt{j}y), \eea and the relevant normalizing constants
\bea \label{co1}&&d_{m,\ell}(m\!-\!1) \! = \! \int_{\bbr^m}\!\!
\frac{(1 \! - \! e^{iy_1})^\ell}{|y|^{2m-1}} dy
\quad \text{\rm ($y_1$ is the first coordinate of $y $)}, \\
\label{co2}&&  \tilde d_{m,k}(m\!-\!1)=\frac{2^{1-m}\pi^{m/2}}{\Gam
(m-1/2)}\intl_0^\infty \frac{(1-e^{-t})^k}{t^{(m+1)/2}}\, dt. \eea

The parameters $k$ and $\ell$ in these formulas will be chosen as
follows. We assume $k$ to be any integer greater than $(m\!-\!1)/2$.
Furthermore, we set $\ell=m-1$ if $m$ is even, and fix any
$\ell>m-1$ if $m$ is odd. Integrals (\ref{co1}) and (\ref{co2}) can
be explicitly evaluated and the following statement holds; see
 [Ru1, pp. 238, 239] and  [SKM, Section
26] for $I^\a f$ with arbitrary $\a$.
\begin{theorem} \label{441} Let $g=I^{m-1} f, \; f \in L^p(\bbr^m), \; 1 \le p <m/(m-1)$.
Then \bea f(x)&=& \frac{1}{d_{m,\ell}(m-1)}\int_{\bbr^m}
\frac{(\Del^\ell_y g)(x)}{|y|^{2m-1}} \,dy  \nonumber\\&=&
\frac{1}{\tilde d_{m,k}(m-1)}\int_{\bbr^m}  \frac{(\tilde\Del^k_y
g)(x)}{|y|^{2m-1}}\, dy, \nonumber\eea
 where $\int_{\bbr^m}=\lim\limits_{\e \to
0}\int_{|y|>\e}$. This limit exists in the $L^p$-norm and in the
a.e. sense. For $f\in C_0 \cap L^p$,  it exists in the $\sup$-norm.
\end{theorem}
 Theorem \ref{441} and (\ref {33a}) give the following result.
\begin{theorem} \label {441a} Let $ f \in L^p(\bbr^m), \; 1 \le p
<m/(m-1)$. We set \be\label{ge} g=(2\pi)^{1-m} \,\stackrel{*}{R_T}
\tilde \vp, \qquad \tilde \vp (a,b)\!=\!\sqrt{1\!+\!|a|^2}\,(R_T f)
(a,b).\ee Then $f$ can be reconstructed by the formulas
\be\label{for1} f(x) \! = \!
\frac{1}{d_{m,\ell}(m\!-\!1)}\int_{\bbr^m} \!\! \frac{(\Del^\ell_y
g)(x)}{|y|^{2m-1}} dy \ee or \be\label{for2} f(x) \!
=\frac{1}{\tilde d_{m,k}(m\!-\!1)}\int_{\bbr^m}  \!\!
\frac{(\tilde\Del^k_y g)(x)}{|y|^{2m-1}} dy, \ee where convergence
of hypersingular integrals is understood as in
 Theorem \ref{441}.
\end{theorem}

\begin{remark} \label{rem1} Theorem \ref{441a} includes continuous functions, the
behavior of which at infinity is specified accordingly.  We denote
\be C_\del(\bbr^m)=\{f \in C(\bbr^m): \; f(x)=O(|x|^{-\del}) \},
\quad \del>0. \nonumber\ee If $\del>m-1$, then there exists $p$ such
that $m/\del<p<m/(m-1)$. Hence  inversion formulas (\ref{for1}) and
(\ref{for2}) are applicable to any $f\in C_\del(\bbr^m)$ with
$\del>m-1$.
\end{remark}

\subsubsection{ Powers of ``minus-Laplacian''} Another series of inversion
formulas can be obtained using integer powers of the operator
$-\Del$, where $\Del$ is the Laplace operator on $\bbr^m$.
\begin{definition} For $\lam \in (0,1),$ let ${\rm Lip}^{loc}_\lam (\bbr^m)$ be
the space of functions $f$ on $\bbr^m$ having the following
property: for each bounded domain $\Om $ in $ \bbr^m$, there is a
constant $A>0$ such that \be |f(x)-f(y)| \le A|x-y|^\lam \qquad
\forall x,y \in \bar{\Om} \quad \text{\rm (the closure of $\Om$)}.
\nonumber\ee We denote \be C^*_\del(\bbr^m)\!=\!\{f: \;  f \!\in
\!C_\del (\bbr^m)\cap {\rm Lip}^{loc}_\lam (\bbr^m)\quad \text{\rm
for some $\lam \in (0,1)$} \}. \nonumber\ee
\end{definition}
\begin{theorem}  \label{invek2} Let $$g=(2\pi)^{1-m} \,\stackrel{*}{R_T} \tilde \vp,
\qquad \tilde \vp (a,b)\!=\!\sqrt{1\!+\!|a|^2}\,(R_T f) (a,b).$$

 {\rm (i)} If $m$ is odd, $\del>m-1$, and $f \in C^*_\del(\bbr^m),$
 then
 \be \label{f1}f(x)= (-\Del)^{(m-1)/2}\,g(x).\ee

{\rm (ii)} If $m$ is even,  $\del>m-1$,  and $f \in C_\del(\bbr^m)$,
 then
  \be \label{f2} f(x)=\frac{2}{\sig_m}\int_{\bbr^m}
\frac{(-\Del)^{m/2-1}g (x)-(-\Del)^{{m/2-1}}g (x-y)}{|y|^{m+1}} \,
dy, \ee where $\int_{\bbr^m}=\lim\limits_{\e \to 0}\int_{|y|>\e}$
uniformly in $x \in \bbr^m$.

Moreover,  if $ m\ge 4, \; f \in C^*_\del(\bbr^m), \; \del>m-1$,
then \be \label{f3} f(x)=-(-\Del)^{(m-2)/2}(I^1 \Del g) (x).\ee

All derivatives in (\ref {f1})-(\ref {f3}) exist in the classical
 sense.
\end{theorem}
\begin{proof}
 These statements are consequences of known facts for potentials
and  singular integrals. In  our case $g=I^{m-1} f$; see
(\ref{33a}).

(i) To ``localize'' the problem, let $x \in B_R=\{x: \; |x|<R \}$
and choose $\chi (x) \in C^\infty (\bbr^m)$ so that \[ 0 \le \chi
(x) \le 1, \; \chi (x)\equiv 0 \; \text{\rm if} \; |x| \le R+1, \;
\text{\rm and} \; \chi (x) \equiv 1 \; \text{\rm if} \; |x|\ge
R+2.\] We have $f=f_1 +f_2, \; f_1=\chi f, \; f_2=(1-\chi)f$, \be
f_1 (x) \! = \! \left\{ \!
 \begin{array} {ll} \!  0 \!  & \mbox{if $ |x| \!  \le \!  R \! + \! 1,$}\\
  \!  f(x) \!  & \mbox{if $  |x| \! \ge  \! R \! + \! 2,$}\\
  \end{array}
\right. \quad  f_2 (x) \! = \!  \left\{ \!
 \begin{array} {ll}  \! f(x) \!  & \mbox{if $  |x| \!  \le R \! + \! 1,$}\\
  \!  0  \! & \mbox{if $  |x| \! \ge \!  R \! + \! 2.$}\\
  \end{array}
\right. \nonumber\ee Let $g=g_1 + g_2, \; g_1=I^{m-1}f_1, \;
g_2=I^{m-1}f_2$. Then
 $g_1 \in C^\infty (B_R)$, and for all multi-indices
$\gam$,
\[
\partial^\gam g_1 (x)=\frac{1}{\gam_m (m-1)} \intl_{|y|>R+1} f_1 (y)
\,
\partial^\gam |x-y|^{-1} dy. \] In particular, for $m$ odd, we
get $(-\Del)^{(m-1)/2} g_1 (x) =0$. The function $g_2$ belongs at
least to $C^{m-2} (B_R)$, and differentiation is possible under the
sign of integration; see, e.g., \cite [Section 1(6)]{Vl}.

Hence, for $m$ odd, $(-\Del)^{(m-3)/2}g_2=I^2 f_2$ (the Newtonian
potential over a bounded domain), and (i) follows by Theorem 11.6.3
from \cite [p. 231]{Mi}.

(ii) Consider the case $m$ even. By reasoning from above, \be
(-\Del)^{m/2-1} g(x)=(I^1 f)(x), \nonumber\ee and (\ref{f2}) holds
owing to Remark \ref{rem1}. If $m\ge 4$ then, as in (i), we have
$-\Del g=I^{m-3} f$. Hence $-I^1 \Del g=I^{m-2}f$, and (\ref{f3})
follows.
\end{proof}

\subsubsection{ Wavelet transforms} Riesz
fractional derivatives $\bbd^\a$ can be represented using continuous
wavelet transforms; see \cite[Section 17]{Ru1}. We apply this fact
in the context of the Radon inversion procedure.

\begin{definition} A measure $\mu$ on $\bbr^m$,  $m \ge 2$, is said to
be radial if
$$
\int_{\bbr^m} g(\gam x)\, d\mu (x) = \int_{\bbr^m} g(x) \,d\mu (x)$$
for all rotations $\gam \in SO (n)$ and all $\mu$-integrable
functions $g$. \end{definition}

Let $\mu$ be a radial finite Borel measure on $\bbr^m$ and let \be
 (Wg)(x,t)\equiv (g \ast \mu_t) (x) = \int_{\bbr^m} g(x-ty)\,
d\mu (y)\nonumber\ee be the convolution of  $g$ with the dilated
version of $\mu$. If $\mu$ has a certain number of vanishing moments
and
 obeys some regularity conditions,
 then  $(Wg)(x,t)$ is
 called the {\it continuous wavelet transform of g, generated by the wavelet
 measure $\mu$}. In particular, if $\mu$ is absolutely continuous,
 that is, $d\mu(y)=w (|y|)\, dy$, then
$$(Wg)(x,t)=\frac{1}{t^m}\int_{\bbr^m} \,g(y) w\left
(\frac{|x-y|}{t}\right )\, dy.$$ Formally, \be (\bbd^\a
g)(x)={1\over d_{\mu} (\a)}\int_0^\infty
\frac{(Wg)(x,t)}{t^{1+\a}}\, dt,\qquad d_{\mu}
(\a)=\const,\nonumber\ee that can be easily checked using the
Fourier transform. Below we apply Theorem 17.10 from \cite{Ru1} to
our case, when $\a=m-1$ and $g$ has the form (\ref{ge}). Recall that
by (\ref{33a}), $g=I^{m-1} f$, $f\in L^p(\bbr^m)$.

To state the  result, we assume that $\mu$ is a radial Borel
measure, satisfying \bea \label{cd1} &&\int_{|x| > 1} |x|^\b
\,d|\mu|(x) <
\infty \quad \text {for some}  \; \beta > m-1,  \\
\label{cd2} &&\int_{\bbr^m} x^j \,d\mu (x) = 0 \quad \text {for}\;
|j| = 0, 2, 4, \ \ldots \ ,2[(m-1)/2]. \eea We denote \bea
\label{co3}
 \quad \qquad d_{\mu,m}\! \!\!&=&\! \!\!\frac {\pi^{m/2} 2^{1-m} }
{\sig_{m-1} \Gamma (m-1/2)}\\
&\times& \left\{ \!\!
 \begin{array} {ll}  \displaystyle{\Gamma ((1-m)/2) \int_{\bbr^m} \!|y|^{m-1}
\,d\mu (y)} &\mbox{if $m$ is even,}\\
{}\\ \!\displaystyle{\frac {2 (-1)^{(m+1)/2}}{((m\!-\!1)/2)!}
\,\int_{\bbr^m}\!
|y|^{m-1} \,\log |y|\, d\mu (y)} &\mbox{if $m$ is odd.}\\
  \end{array}
\right.\nonumber\eea

\begin{theorem}  \label{441ab} Let $ f \in L^p(\bbr^m), \; 1 \le p
<m/(m-1)$,$$ g=(2\pi)^{1-m} \,\stackrel{*}{R_T} \tilde \vp, \qquad
\tilde \vp (a,b)\!=\!\sqrt{1\!+\!|a|^2}\,(R_T f) (a,b).$$ If
$d_{\mu,m} \neq 0$, then $f$ can be reconstructed by the formula \be
f(x)\!=\!{1\over d_{\mu,m}}\int_0^\infty \frac{(Wg)(x,t)}{t^m}\,
dt\!=\!\lim\limits_{\e \to 0}{1\over d_{\mu,m}}\int_\e^\infty
\frac{(Wg)(x,t)}{t^m}\, dt,\nonumber\ee where the
 limit exists in the $L^p$-norm and in the a.e. sense. For $f\in
C_0 \cap L^p$,  it exists in the $\sup$-norm.
\end{theorem}

Examples of wavelet measures can be found in \cite[Section
17.4]{Ru1}.

\subsection{The second approach: Ridgelet-like
transforms and the convolution-backprojection method}

\subsubsection{Preliminary discussion} Our  aim is to give precise meaning to the formula
\be\label{opta2}f=(2\pi)^{1-m} \, \stackrel{*}{R_T}{}\!\!^{1-m}\vp,
\qquad \vp=R_T f,\ee where $f\in L^p(\bbr^m)$ and
$\stackrel{*}{R_P}{}\!\!^{1-m}\vp$ is obtained by formal
substitution $\a=1-m$ in the fractional integral in (\ref{dt1}). Of
course, we cannot replace $\a$ by $1-m$  directly. To overcome
this difficulty, we proceed as follows. Recall that
 \be
\frac{|a\cdot x'+b -x_m |}{\sqrt{1+|a|^2}}=\dist (x,h) \nonumber\ee
is the Euclidean distance between the point $x\in \bbr^m$ and the
hyperplane $h=\{y\in \bbr^m: a\cdot y' +b=y_m\}$. Thus,
\be\label{dtax} (\stackrel{*}{R_T}{}\!\!^\a
\vp)(x)=\frac{1}{\gamma_{1}(\a)} \int_{\tilde \bbr^m} \vp(a,b)\,
[\dist (x,h)]^{\a-1} \, \tilde d a\, db.\nonumber\ee The kernel of
this integral operator has a singularity on the hyperplane $h$. Let
us write this kernel  in the form
$$
[\dist (x,h)]^{\a-1}= c_{\a, w} \int^\infty_0 w \left({\dist
(x,h)\over t}\right) \; {dt\over t^{2-\a}}, \; \qquad c^{-1}_{\a, w}
= \intl^\infty_0 w(s)\;
 {ds\over s^{\a}}.
$$
Here  $w(s)$ is a measurable function (it will be specified later)
such that
 $\int^\infty_0 |w(s)|/s^{\a}\,ds<\infty$ and $c^{-1}_{\a, w} \neq 0$.
Changing the order of integration, we obtain \be\label{dtay}
(\stackrel{*}{R_T}{}\!\!^\a \vp)(x)=\frac{c_{\a, w}}{\gamma_{1}(\a)}
\int^\infty_0 \, {(\tilde W \vp)(x, t)\over t^{2-m-\a}}\,
 dt,\ee
 where
\be \label{wavt}(\tilde W \vp)(x,t)=\frac{1}{t^m}\,\int_{\tilde
\bbr^m} \vp(a,b)\,w\left (\frac{\dist (x,h) }{t}\right )\,\tilde d
a\, db.\ee If $w$ is a wavelet function with a certain number of
vanishing moments and a suitable decay at infinity, then
(\ref{dtay}) can be meaningful also for negative $\a$, in
particular, for $\a=1-m$. In this case
 \be\label {ftu}
(\stackrel{*}{R_T}{}\!\!^{1-m}\vp)(x)= c_{m, w}\,\int^\infty_0 \,
{(\tilde W \vp)(x, t)\over t}\,
 dt,\ee
where  $c_{m, w}$ is a suitable constant and convergence of the
integral  is  interpreted in a proper way.

An idea of this approach in a slightly different form was suggested
in \cite [Section 10.7]{Ru1} as a
 generalization of Marchaud's method in fractional calculus.

Note that $f$ can be reconstructed directly from $(\tilde W
R_Tf)(x,t)$ if
 this composition  represents an approximate identity.
This observation is a core of the convolution-backprojection method
 in tomography \cite{Na, Ru4} and  agrees  with the
method of approximative inverse operators in fractional calculus;
cf. \cite{NS, Sam3}.

 Our first task in this section is to show that $w$ can be chosen so that
\be\label{lim} \lim\limits_{t\to 0}\,(\tilde W
R_Tf)(x,t)=c_w\,f(x),\ee where $c_w\neq 0$ is a constant depending
on $w$ and the limit is
 understood in the $L^p$-norm and in the
almost everywhere sense. Then we replace the passage to the limit in (\ref{lim})  by integration against the dilation-invariant measure
$dt/t$ and get \be\label {saim} \int_0^\infty
 \frac{(\tilde W R_Tf)(x,t)}{t}\, dt=c'_w\,f(x), \qquad c'_w=\const\neq 0.\ee  Formula (\ref{saim})
agrees with equality (\ref{opta2}), in which  $\stackrel{*}{R_T}{}\!\!^{1-m}\vp$
is represented by (\ref{ftu}).

The convolution (\ref{wavt})
 has the same nature as the ridgelet transform in \cite{C1, C2}.
 A detailed exposition of this topic for totally geodesic Radon transforms
  can be found in \cite{Ru8}.

\subsubsection{Justification of (\ref{lim})}

\begin{lemma} \label{stei} Let $f \in L^p(\bbr^m), \;  1 \le p < m/(m-1), \; m\ge 2$.
Suppose that \be \label {beha} w(s)\,=\, \left\{
\begin{array}{ll} O(s^{\vare_1-1})  &  \mbox{\rm if $s \le 1,$ } \\
 O(s^{-m/p'-\vare_2}) &  \mbox{\rm if $s\ge 1, \quad 1/p+1/p'=1,$ }
\end{array}
\right. \ee for sufficiently small $\vare_1, \vare_2 >0$ and let
\be\label{kern} k(y)=\frac{\sig_{m-2}}{|y|^{m-2}}\int_0^{|y|} (|y|^2
-s^2)^{(m-3)/2}w(s)\,ds,\ee  $k_t(y)=t^{-m}k(y/t), \;t>0$. Then
\be\label{convo} (\tilde W R_Tf)(x,t)=(f\ast k_t)(x)\ee for almost
all $x$.
\end{lemma}
\begin{proof} Changing variables and using Fubini's theorem, we
obtain
$$
(\tilde W R_Tf)(x,t)=\int_{\bbr^m} f(x-ty)\, k(y)\, dy, $$ $$
k(y)=\int_{\bbr^{m-1}} w \left (\frac{|y\cdot
(a+e_m)|}{\sqrt{1+|a|^2}} \right )\, \tilde d a.
$$
By (\ref{teq2}),$$
 k(y)=\int_{S^{m-1}_+} w(|y\cdot\theta|)\, d\theta=\frac{1}{2}\int_{S^{m-1}} w(|y\cdot\theta|)\,
 d\theta,$$ and (\ref{hva37}) yields
$$
 k(y)=\sig_{m-2}\int_0^1 (1-t^2)^{(m-3)/2}w(|y|t)\,dt.$$
 This gives (\ref{kern})-(\ref{convo}). To complete the proof we
 must justify application of Fubini's theorem. It suffices to show
 that $(\tilde W R_Tf)(x,t)$ is finite a.e. when $f$ and $w$ are replaced by
 their absolute values. We split the integral
\be\label{spli} I(x,t)=\frac{1}{t^m}\,\int_{\tilde \bbr^m} (R_T
|f|)(a,b)\,\left |w\left (\frac{\dist (x,h) }{t}\right )\right
|\,\tilde d a\, db \ee into two pieces $I_1 +I_2$, where in $I_1$ we
integrate over the set $\{
 (a,b): \dist (x,h)/t <1\}$ and in $I_2$  over the set $\{
 (a,b): \dist (x,h)/t >1\}$. Owing to (\ref{beha}) and (\ref{dt1}),
$$
 I_i\le c \,t^{{1-\b_i -m}}(\stackrel{*}{R_T}{}\!\!^{\b_i} R_T
|f|)(x), \qquad i = 1, 2;$$ $$ \beta_1=\vare_1, \qquad
\beta_2=1-m/p' -\vare_2, \qquad c = \const.$$
 By
Theorem \ref{t31}, $I_i$ are dominated by multiples of Riesz
potentials, namely, \be\label {riep} I_i\le c \,t^{1-\b_i -m}
(I^{\beta_i+m-1} |f|)(x).\ee The latter are finite for almost all
$x$ (see, e.g., \cite[Chapter V, Section 1.2]{St1}) because, for
sufficiently small $\vare_1$
 and $\vare_2 >0$ we have  $0<\beta_i+m-1 <m/p$. This completes the
 proof.
\end{proof}
Owing to Lemma \ref {stei}, the standard machinary of approximation
to the identity (see,
 \cite[Chapter III, Section 2.2]{St1})) implies the following
 statement.

\begin{theorem} \label{cbp} Let $ 1 \le p < m/(m-1), \; m\ge 2$.
Suppose that $w(s)$ satisfies (\ref{beha})  and the corresponding
  kernel $$ k(y)=\frac{\sig_{m-2}}{|y|^{m-2}}\int_0^{|y|} (|y|^2
-s^2)^{(m-3)/2}w(s)\,ds$$  has a radial decreasing majorant in
 $L^1(\bbr^m)$. If $f \in L^p(\bbr^m)$,
then \be\label{510} \lim\limits_{t\to 0}\frac{1}{t^m}\,\int_{\tilde
\bbr^m} (R_Tf)(a,b)\,w\left (\frac{|a\cdot x'+b -x_m
|}{t\sqrt{1+|a|^2}}\right )\,\tilde d a\, db= \gamma\,f(x),\ee
\be\label{gaam} \gamma= \int_{\bbr^m} k(y)\, dy,\ee where the limit
 exists in the $L^p$-norm and in the a.e. sense. If, moreover,
$f \in C_0 (\bbr^m)\cap L^p(\bbr^m)$, then the convergence in
(\ref{510}) is uniform on $\bbr^m$.
\end{theorem}

To give an example of a function $w$ satisfying Theorem \ref{cbp},
we invoke Riemann-Liouville fractional integrals (\ref{ril})
 and write $k(y)$ in the form \be
k(y)=\pi^{(m-1)/2}r^{2-m}(I^{(m-1)/2}_{0+} [s^{-1/2}
w(\sqrt{s})])(r^2), \qquad r=|y|.\nonumber\ee Then
 \be \gamma\equiv \int_{\bbr^m}
k(y)\, dy=c\,\int^\infty_0 \lambda (t)\,dt, \nonumber\ee $$ \lambda
(t)=(I^{(m-1)/2}_{0+} [s^{-1/2} w(\sqrt{s})])(t), \quad
c=\frac{\pi^{m-1/2}}{\Gam (m/2)}.$$

\begin{example} \label{exa1}(cf. [SS], p. 338).  Let
\be
\kappa_{\a,\ell}(s)=\big(\dfrac{d}{ds}\big)^\ell\;\dfrac{s^\ell}{(s+i)^{1+\a}},
\qquad \ell \in \bbn,\quad \a > 0,\nonumber\ee $$
  \lam_{\a,\ell} (t) \equiv
(I^{\a}_{0+}\kappa_{\a,\ell})(t)=\frac{i^{\ell-\a} \, \ell!}{\Gam
(1+\a)} \,\frac {t^\a}{(t+i)^{\ell+1}}.
$$
Both functions are  integrable on $(0, \infty)$ for $\ell>\a$ and
\be\label {iny}\int^\infty_0\lambda_{\a, \ell} (t)\, dt = \Gamma(\ell-\a).\ee We
set $\a = (m-1)/2$. Then for any $\ell> (m-1)/2$ the function $$
 w (s)=s\kappa_{(m-1)/2, \ell}(s^2)$$ can be used in Theorem
 \ref{cbp}.
The corresponding constant $\gamma$ in (\ref{gaam}) will be computed
as
$$\gamma=c\,\int^\infty_0 \lambda (t)\,dt=
 c\,\int^\infty_0 \lambda_{(m-1)/2, \ell} (t)\,dt=
 \frac{\pi^{m-1/2}\, \Gamma(\ell-(m-1)/2)}{\Gam (m/2)}.$$
\end{example}

\subsubsection{Justification of (\ref{saim})} Let us
 represent the
corresponding truncated integral as an approximate identity.
\begin{lemma} \label{stei2} Let $f \in L^p(\bbr^m), \;  1 \le p < m/(m-1), \; m\ge 2$.
Suppose that $w$ satisfies (\ref{beha})
 and let
\be\label{kern1} g(y)=\frac{\sig_{m-2}}{(m-1) |y|^{m}}\int_0^{|y|}
(|y|^2 -s^2)^{(m-1)/2}w(s)\,ds,\ee  $g_\e(y)=\e^{-m}g(y/\e),
\;\e>0$. Then \be\label{convo2} \int_\e^\infty
 \frac{(\tilde W R_Tf)(x,t)}{t}\, dt
=(f\ast g_\e)(x)\ee for almost all $x$.
\end{lemma}
\begin{proof} The result follows  from (\ref{convo}) if we change variables and apply Fubini's
theorem. To justify application of Fubini's theorem, we note that
the left-hand side of (\ref{convo2}) is  finite a.e. when $f$ and
$w$ are replaced by $|f|$ and $|w|$, respectively. Indeed, let
$$
\tilde I_{\e}(x)=\int_\e^\infty
 \frac{I(x,t)}{t}\, dt,$$ where $ I(x,t)$ is defined by
 (\ref{spli}). By (\ref{riep}),
 $$
\tilde I_{\e}(x)\le c\,\sum\limits_{i=1}^2 (I^{\beta_i+m-1}
|f|)(x)\,\int_\e^\infty\frac{dt}{t^{\b_i +m}},$$
 which is finite a.e. because $\b_i
+m>1$.

It remains to perform simple calculations. We have $$ \int_\e^\infty
 \frac{(\tilde W R_Tf)(x,t)}{t}\, dt=\int_\e^\infty
 \frac{(f\ast k_t)(x)}{t}\, dt=\left (f\ast \int_\e^\infty
 \frac{k_t}{t}\, dt\right )(x).$$
Set $k(y)= k_0 (|y|)$. Then, by (\ref{kern}),
$$
\int_\e^\infty
 \frac{k_t (y)}{t}\, dt=\int_\e^\infty k_0 \left (\frac{|y|}{t}\right )\,
 \frac{dt}{t^{m+1}}=\frac{1}{\e^m}\, g\left ( \frac{y}{\e} \right
 ),$$ where $g$ has the form (\ref{kern1}).
\end{proof}

Lemma \ref  {stei2} implies the following inversion result.
\begin{theorem} \label{cbpx} Let  $f \in L^p(\bbr^m)$, $ \;1 \le p < m/(m-1), \; m\ge 2$.
Suppose that $w(s)$ satisfies (\ref{beha}) and the corresponding
  kernel $$g(y)=\frac{\sig_{m-2}}{(m-1) |y|^{m}}\int_0^{|y|}
(|y|^2 -s^2)^{(m-1)/2}w(s)\,ds$$  has a radial decreasing majorant in
 $L^1(\bbr^m)$. If
 \be\label{waa} (\tilde W \vp)(x,t)=\frac{1}{t^m}\,\int_{\tilde
\bbr^m}\vp(a,b)\,w\left (\frac{|a\cdot x'+b -x_m
|}{t\sqrt{1+|a|^2}}\right )\,\tilde d a\, db,\nonumber\ee
 where $\vp=R_T f$, then \be\label{510x} \lim\limits_{\e\to 0}\int_\e^\infty \frac{(\tilde W \vp)(x,t)}{t}\, dt
= \gamma_1\,f(x),\ee \be\label{gaamx} \gamma_1= \int_{\bbr^m} g(y)\,
dy,\nonumber\ee  the limit
 being understood in the $L^p$-norm and in the a.e. sense. If, moreover,
$f \in C_0 (\bbr^m)\cap L^p(\bbr^m)$, then the convergence in
(\ref{510x}) is uniform on $\bbr^m$.
\end{theorem}

To give an example of a function $w$ satisfying Theorem \ref{cbpx},
we  write $g(y)$
 in the form \be
g(y)=\frac{\pi^{(m-1)/2}}{2r^{m}}\,(I^{(m+1)/2}_{0+} [s^{-1/2}
w(\sqrt{s})])(r^2), \qquad r=|y|.\nonumber\ee Then
 \be \label {ppo}\gamma_1\equiv \int_{\bbr^m}
g(y)\, dy=c\,\int^\infty_0 \lambda_1 (t)\,dt, \quad
c=\frac{\pi^{m-1/2}}{2\Gam (m/2)},\ee \be\label{kouo} \lambda_1
(t)=t^{-1}(I^{(m+1)/2}_{0+} [s^{-1/2} w(\sqrt{s})])(t).\ee

\begin{example} \label{exa2}  Let
$$
 w (s)=s\kappa_{(m+1)/2, \ell}(s^2), \qquad \ell> (m-1)/2;$$
see Example \ref {exa1}. Then
\bea
\lambda_1 (t)&=&t^{-1}\lam_{(m+1)/2,\ell} (t)=\frac{i^{\ell-(m+1)/2} \, \ell!}{\Gam
(1+(m+1)/2)} \,\frac {t^{(m-1)/2}}{(t+i)^{\ell+1}}\nonumber\\
&=&\frac{2}{i(m+1)}\, \lam_{(m-1)/2,\ell} (t)\quad (\in L^1 (0,\infty)).\nonumber\eea
Hence, by (\ref{iny}),
\bea
\gamma_1&=&c\,\int^\infty_0 \lambda_1 (t)\,dt=\frac{2c}{i(m+1)}\,
\int^\infty_0 \lam_{(m-1)/2,\ell} (t)\,dt\nonumber\\
&=&\frac{\pi^{m-1/2}\, \Gam (\ell -(m-1)/2)}{i(m+1)\Gam (m/2)}.\nonumber\eea
\end{example}

\begin{remark} The assumption $g\in L^1(\bbr^m)$ in Theorem \ref{cbpx} can be fulfilled only
if the  generating function $w$ has
vanishing moments. Hence $(\tilde W \vp)(x,t)$ is a
 wavelet-like transform (or a ridgelet-like transform)
 and (\ref{510x}) has the same nature as the classical Calder\'on identity; cf. \cite{FJW}.

 Furthermore, by Lemma \ref {frai}, for any function $w$,
 satisfying
\be \int^\infty_0 s^{2j} w(s) \,ds=0 \; \; \; \forall \; j = 0, 2,
\dots, 2 [(m-1)/2],\nonumber\ee \be \int^\infty_1 s^{\beta} |w(s)|
\,ds<\infty \; \; \hbox{\rm for some } \; \; \beta >
m-1,\nonumber\ee the corresponding function (\ref{kouo}) has a
decreasing integrable majorant.  In this case the constant
$\gamma_1$ in (\ref{ppo}) has the form \be\label{219} \gamma_1=
c\,\left\{ \!\!
 \begin{array} {ll}  \displaystyle{2\Gam \left(\frac{1-m}{2}\right )\,
 \int_0^\infty s^{m-1}w(s)\,ds} &\mbox{ if $m$ is even},\\
{}\\ \!\displaystyle{\frac{4(-1)^{(m+1)/2}}{((m-1)/2)!}\int_0^\infty
s^{m-1}w(s)\,\log s\, ds}
&\mbox{ if $m$ is odd}.\\
 \end{array}
\right.\nonumber\ee
Such functions $w(s)$,
having fast
 decay as $s \to 0$ and $ \infty$, can be constructed using Example 1.6 from \cite{Ru3}.
\end{remark}

A similar remark addresses to the assumption $k\in L^1(\bbr^m)$ in
Theorem \ref{cbp}.

\section{Transversal Radon transforms on the  Semyanistyi-Lizorkin
spaces}

Our aim in this section is to bring light to Strichartz's inversion
formulas (\ref{Siv}) in the general context  of the transversal
Radon transform and to determine spaces, say, $X$ and $Y$, of smooth
functions such that $R_T$ acts from $X$ onto $Y$ as an isomorphism.

\subsection{Isomorphism of $R_T$} \label{sed}  To start with, we observe that  if $f\in L^1(\bbr^{m})$
is nonnegative, then $\vp =R_Tf$ is integrable on $\tilde\bbr^{m}$
only if $f\equiv 0$. This is obvious from the equality
$$\int_{\tilde\bbr^m} \!\!\vp(a,b)\,da
db\!=\!\int_{\bbr^{m-1}}\!\!\!\!\!da\int_{\bbr}
db\int_{\bbr^{m-1}}\!\!\! f(x', a\cdot
x'+b)\,dx'=\!||f||_1\int_{\bbr^{m-1}}\!\!\!\!da.$$
Thus, the
 relation $R_T f\in L^1(\tilde \bbr^{m})$ is
 possible only if $f$ is sign-changing. It means that embedding $R_T [S(\bbr^{m})] \subset S(\tilde\bbr^{m})$
 does not hold, although, $R_T [S(\bbr^{m})] \subset
 C^\infty(\tilde\bbr^{m})$.

A close  situation is  known in the theory of operators
of the potential type, where
Semyanistyi-Lizorkin spaces of ``good" functions come into play
 and lead to deep results; see \cite{Ru1, Sam3, SKM}. These
spaces were introduced by Semyanistyi \cite{Se1} and essentially generalized by
Lizorkin \cite{Liz1}-\cite{Liz3} and Samko \cite{Sam1}-\cite{Sam3}. They
incorporate Schwartz functions, vanishing on a given set, and their
Fourier images with zero moments.

 Given a function $f (x)\equiv f (x', x_m)$ on $\bbr^m$, let $$
(\F f)(y)=\int_{\bbr^m}f(x) e^{ix\cdot y}dx$$
 be the Fourier transform of  $f$. We denote by $\F_1$ and $\F_2$  the similar transforms in the $x'$-variable and the $x_m$-variable, respectively.
Let $\Psi\equiv \Psi (\bbr^m)$
  be the subspace of functions $\psi \in S(\bbr^m)$
vanishing with all derivatives $\partial_m^k \psi$, $k\in \bbz_+$,
on the hyperplane $x_m=0$.
 We denote by
 $\Phi\equiv \Phi(\bbr^m)=\F[ \Psi (\bbr^m)]$   the
  Fourier image of $ \Psi (\bbr^m)$.  The spaces $ \Psi (\tilde\bbr^m)$ and $ \Phi (\tilde\bbr^m)$ have the same meaning.

The following auxiliary statements are immediate consequences of preceding definitions.
 For convenience of the reader we present them with
proofs.

\begin{proposition} \label{jikl1} ${}$\hfill

\noindent {\rm (i)} A Schwartz function $\phi$ belongs to
$\Phi(\bbr^m)$ if and only if
  \be\label{kssi1} \int_{\bbr} \phi(x', x_m)\, x_m^{k}
 \,dx_m=0 \qquad
 \forall k\in \bbz_+, \quad \forall x' \in\bbr^{m-1}.\ee

\noindent {\rm (ii)} $\F_1$ is an automorphism of $\Psi (\bbr^m)$.

\noindent {\rm (iii)}  $\F_2$ acts as an isomorphism from $\Psi
(\bbr^m)$ onto $\Phi (\bbr^m)$.

\noindent {\rm (iv)} The map $\partial^{\bf j}: \Psi(\bbr^m) \to
\Psi(\bbr^m)$ is continuous for every ${\bf j}\in \bbz_+^m$.

  \end{proposition}
\begin{proof} {\rm (i)} Let $\phi \in \Phi$, that is, $\phi=\F\psi$ for some $\psi \in
\Psi$. Then for any $k\in \bbz_+$, \bea (\partial_m^k \psi)(y',
y_m)=\frac{1}{2\pi} \F_1^{-1}\Big [\int_{\bbr}\phi(x', x_m)\,
(-ix_m)^k e^{-ix_m y_m}dx_m \Big ].\nonumber\eea Owing to
injectivity of $\F_1^{-1}$,
 the result follows.

\noindent {\rm (ii)}  Let $\psi \in \Psi, \; \psi_1=\F_1\psi.$ Then
for any $k\in \bbz_+$,
$$
(\partial_m^k \psi_1)(y', 0)=\int_{\bbr^{m-1}} (\partial_m^k \psi)
(x',0)\, e^{ix'\cdot y'}dx'=0, \quad \mbox {\rm i.e.,} \quad \psi_1
\in \Psi.
$$
Conversely, we have $\psi=\F_1 \psi_2$, where $$\psi_2 (y', y_m)=(\F_1^{-1}
\psi) (y', y_m)= (2\pi)^{1-m}(\F_1 \psi) (-y', y_m)\in \Psi.$$

\noindent {\rm (iii)} If $\psi \in \Psi$, then $\F_2 \psi=\F\psi_0$,
$\psi_0=\F^{-1}\F_2 \psi=\F_1^{-1} \psi\in \Psi$ by {\rm (ii)}.
Hence  $\F_2 \psi \in \Phi$. Conversely, let $\phi \in \Phi$, that
is, $\phi=\F\psi$, $\psi \in \Psi$. Then $\phi=\F_2 \F_1\psi$ with
$\F_1\psi\in \Psi$.

\noindent {\rm (iv)}  Let $\psi \in \Psi$, that is, $\psi
(y)=(\F^{-1} \vp)(y)=(2\pi)^{-m}(\F \phi) (-y)$, $\phi \in \Phi$.
Then $(\partial^{\bf j}\psi)(y)=(2\pi)^{-m}(-1)^{|{\bf
j}|}\,(\partial^{\bf j}\F \phi) (-y)$, and therefore, by (i), \bea
(\partial^{\bf j}\psi)(y',0)&=&(2\pi)^{-m}(-1)^{|{\bf
j}|}\int_{\bbr^{m-1}}(ix')^{\bf j'}\, e^{-ix' \cdot
y'}dx'\nonumber\\&\times& \int_{\bbr}(ix_m)^{j_m}\,\phi (x', x_m)\,
dx_m=0.\nonumber\eea The continuity of the map $\partial^{\bf j}:
\Psi \to \Psi$ follows from its continuity in the topology of
$S(\bbr^m)$.
\end{proof}

\begin{theorem} \label{ish}The transversal Radon transform $R_T$ acts as an isomorphism from $\Phi (\bbr^m)$ onto $\Phi
(\tilde \bbr^m)$. \end{theorem}

The proof of Theorem \ref{ish} relies on the following two lemmas.

\begin{lemma}\label{ish1} For any $f \in L^1 (\bbr^m)$,
\be  \label{ash1}\F_2 [(R_T f) (a, \cdot)](\xi)=(\F f)(-a\xi, \xi),
\qquad \forall (a, \xi)\in \tilde
 \bbr^m.\ee
\end{lemma}
This formula can be easily obtained by direct calculation, using Fubini's theorem. Introducing a ``mixing" map
 \be
(\Lam u)(a, \xi)= u(-a\xi, \xi),\ee we formally have \be\label{roz}
 R_T f=\F_2^{-1}\Lam\F f.\ee

\begin{lemma}\label{ish3} The map $\Lam $ acts as an isomorphism from $\Psi (\bbr^m)$  onto $\Psi (\tilde \bbr^m)$.
\end{lemma}
\begin{proof} STEP 1. Let $u (y)\equiv  u
(y', y_m)\in \Psi (\bbr^m)$. For every $p,q\in\bbz_+$ there is a
constant $c_{p,q}$ such that \be \label{ash2} |u (y)|\le c_{p,q}\,
|y_m|^{2p}\, (1+|y|^2)^{-q}\qquad \forall y\in \bbr^m;\ee cf.
\cite[Ch. 2, Section 4]{Sam3}. If $|y_m|>1$ this inequality is obvious
and holds for every Schwartz function. In the case $|y_m|\le 1$ the
estimate
  can be obtained by making use of  Taylor's formula with integral remainder.

We first show that \be\label{ook}v(a, \xi)\!=\!u(-a\xi, \xi)\!\in
\!\Psi (\tilde \bbr^m).\ee  Clearly, $v\! \in \!C^\infty (\tilde
\bbr^m)$ and a simple calculation shows that every derivative
$\partial ^{\bf j}v$, $\bf j\in \bbz^m_+$, can be represented as a
finite sum of expressions of the form
$$
U_{{\bf i}}(a, \xi)=Q_{{\bf i}}(a, \xi)\,(\partial ^{\bf i}u)(-a\xi,
\xi), \qquad \bf i\in \bbz^m_+,$$   $Q_{{\bf i}}$  being certain
polynomials. Since differentiation preserves the space $\Psi
(\bbr^m)$ (see Proposition \ref{jikl1} (iv)), $U_{{\bf i}}(a,
0)\equiv 0$, and multiplication by a polynomial preserves $\Psi
(\tilde \bbr^m)$, it suffices to show that the function $v_{{\bf
i}}(a, \xi)\!=\!(\partial ^{\bf i}u)(-a\xi, \xi)$ is rapidly
decreasing. In other words, we have to show that for any $k \in
\bbz_+$ there is a constant $c_k$ such that \be\label {knb}
(1+|a|^2+\xi^2)^k |v_{{\bf i}}(a, \xi)| \le c_k \quad \forall (a,
\xi) \in \tilde \bbr^m. \ee This would imply  $v\! \in \!S (\tilde
\bbr^m)$. We have
$$l.h.s.=\sum\limits_{i=0}^{k}{k \choose i}V_{k,i}, \quad
V_{k,i}=(1+|a|^2)^i \,(\xi^2)^{k-i}\,|v_{{\bf i}}(a, \xi)|.$$ Owing
to
 (\ref{ash2}) with $p=k$ and $q=2k-i$, we obtain
 $$
V_{k,i}\le c_{k,2k-i}\, \Big
(\frac{\xi^2+|a\xi|^2}{1+\xi^2+|a\xi|^2}\Big )^i \, \Big
(\frac{\xi^2}{1+\xi^2+|a\xi|^2}\Big )^{2k-2i}< c_{k,2k-i}.
$$
This gives (\ref{knb}) with $c_k=\sum\limits_{i=0}^{k}{k \choose
i}c_{k,2k-i}$, which implies (\ref{ook}).

 STEP 2. Now our task is to prove that every function $v\in\Psi (\tilde \bbr^m)$ has the form
 $v(a, \xi)\!=\!u(-a\xi, \xi)$ for some function $u\in \Psi
 (\bbr^m)$. We set
 \be\label {kiu}
  u (y)\equiv  u
(y', y_m)=v(-y'/y_m, y_m), \quad y_m\neq 0; \qquad u (y',0)=0.\ee
 Let us show that this function belongs to $\Psi
 (\bbr^m)$.
As above, for every $p,q\in\bbz_+$ there is a constant $c_{p,q}$
such that \be \label{ash4} |v(a, \xi))|\le c_{p,q}\, |\xi|^{2p}\,
(1+|a|^2+\xi^2)^{-q}\qquad \forall (a, \xi)\in \tilde\bbr^m.\ee
Furthermore, for $y_m\neq 0$, every derivative $(\partial ^{\bf
j}u)(y', y_m)$, $\bf j\in \bbz^n_+$, is a finite sum of expressions
of the form
$$\tilde U_{{\bf i}}(y', y_m)=\tilde Q_{{\bf i}}(y', 1/y_m)\,(\partial ^{\bf i}v)(-y'/y_m, y_m),
\qquad \bf i\in \bbz^m_+,$$ where  $\tilde Q_{{\bf i}}$ are certain
polynomials. Again, since differentiation preserves the space $\Psi
(\tilde \bbr^m)$ and multiplication by $\tilde Q_{{\bf i}}(y',
1/y_m)$ preserves $\Psi (\bbr^m)$, it suffices to show that the
function $u_{{\bf i}}(y)\!=\!(\partial ^{\bf i}v)(-y'/y_m, y_m)$ is
rapidly decreasing when $|y|\to \infty$ and $y_m \to 0$.

Let us check that for any $r,s \in \bbz_+$ there is a constant
$c_{r,s}$ such that \be\label {knb2} y_m^{-2s}(1+|y|^2)^r |u_{{\bf
i}}(y)| \le c_{r,s} \quad \forall y \in \bbr^m. \ee We have
$l.h.s.=\sum\limits_{k=0}^{r}{r \choose k}U_{r,k}$, where
$$
U_{r,k}=y_m^{-2s}\,(1+y_m^2)^k \,(|y'|^2)^{r-k}\,|(\partial ^{\bf
i}v)(-y'/y_m, y_m)|.$$ Owing to
 (\ref{ash4}) with $p=s$ and $q=2r-k$,
 \bea
U_{r,k}&\le& c_{s,2r-k}\, \Big
(\frac{1+y_m^2}{1+y_m^2+|y'|^2/y_m^2}\Big )^k \nonumber\\
&\times& \Big (\frac{|y'|^2/y_m^2}{1+y_m^2+|y'|^2/y_m^2}\Big )^{r-k}
 \, \Big (\frac {y_m^2}{1+y_m^2+|y'|^2/y_m^2}\Big )^{r-k}
< c_{s,2r-k}. \nonumber\eea This gives (\ref{knb2}) with
$c_k=\sum\limits_{k=0}^{r}{r \choose k} c_{s,2r-k}$, which implies
$u\in \Psi
 (\bbr^m)$.

 To complete the proof we note that the equality $v(a, \xi)\!=\!u(-a\xi, \xi)=(\Lam u)(a, \xi)$ for the function
 $u$ defined by (\ref{kiu}) is obvious. Thus the map
$\Lam : \, \Psi (\bbr^m)\to \Psi (\tilde \bbr^m)$ is surjective and
the inverse map $$\Lam^{-1} : \,\Psi (\tilde \bbr^m)\to \Psi
(\bbr^m)$$  is well-defined by \be\label {kiu1}
  (\Lam^{-1}v) (y)\!=\!v(-y'/y_m, y_m), \; y_m\!\neq \!0; \qquad (\Lam^{-1}v)(y',0)\!=\!0.\ee
Both maps $\Lam$ and $\Lam^{-1}$ are continuous in the topology of
the Schwartz space. This can be proved directly using (\ref{ash2})
and (\ref{ash4}).
\end{proof}

\noindent{\it Proof of Theorem \ref{ish}}.  Since the maps
$$\F: \Phi (\bbr^m)\to \Psi (\bbr^m), \quad \Lam :  \Psi (\bbr^m)\to
\Psi (\tilde \bbr^m), \quad \F_2^{-1}: \Psi (\tilde \bbr^m)\to \Phi
(\tilde \bbr^m)$$ are isomorphisms, then formula (\ref{roz}) is
well-justified on functions $f \in \Phi (\bbr^m)$, and the result
follows.

\subsection{Inversion formulas}

A pointwise inversion formula \be\label{roz1}
 f=\F^{-1}\Lam^{-1}\F_2 R_T f,\qquad f\in \Phi (\bbr^m),\ee
in terms of the Fourier transforms follows from (\ref{roz}). Below
we obtain alternative inversion formulas, that do not contain the
Fourier transform. To this end we invoke
 the  partial Riesz potential  \be \label{15p}(I_2^\a f)(x)=\frac{1}{\gamma_1(\a)}
\intl_{\bbr}
 \frac{f(x',y_m)\,dy_m}{|x_m-y_m|^{1-\a}},\quad
  \gamma_{1}(\a)=
  \frac{2^\a\pi^{1/2}\Gamma(\a/2)}{\Gamma((1-\a)/2)},
  \ee
$$Re \, \a >0, \qquad \a\neq 1,3,5,\ldots \,.$$ This operator
is an automorphism of $\Phi(\bbr^m)$ \cite{SKM} and \be\F[I_2^\a
f](y)=|y_m|^{-\a} \F[f](y), \qquad f \in \Phi.\ee  The last relation
extends $I_2^\a f$, $f \in \Phi$, to all $\a \in \bbc$ as an entire
function of $\a$. We also introduce a backprojection operator
$\tilde R_T$\footnote{We do not call it the {\it dual operator}
because the latter has  another meaning in our paper.}, which
sends functions on  $\tilde \bbr^m$ to functions on
$\bbr^m$  by the formula \be\label{bpro} (\tilde R_T
g)(x)=\int_{\bbr^{m-1}} g(a, -x'\cdot a+ x_m)\,da=(R_T g)(-x',
x_m).\ee

\begin{theorem}\label {pzp} If $f \in
\Phi(\bbr^m)$, then for  any complex $\a$ and $\b$, \be\label{cab}
(I_2^\a \tilde R_T I_2^\b R_T  f)(x)=(2\pi)^{m-1} \,(I_2^{\a+\b+m-1}
f)(x), \qquad x \in \bbr^m.\ee \end{theorem}
\begin{proof} Let $g=I_2^\b R_T  f$. Then $g \in \Phi(\tilde\bbr^m)$ and
(\ref{ash1}) yields \bea (\F_2 I_2^\a [(\tilde R_T g)(x',
\cdot)])(\eta)&=& |\eta|^{-\a}(\F_2 [(\tilde R_T g)(x',
\cdot)])(\eta) \nonumber\\ \label{ash7}&=&|\eta|^{-\a}(\F g)(x'\eta,
\eta). \eea Furthermore, denoting by $\F_1$ the Fourier transform in
the first $m-1$ variables, we have \bea (\F g)(x'\eta, \eta)&=&
(\F_1 \{\F_2 [g(a, \cdot)] (\eta) \})
(x'\eta)\nonumber\\
&=& (\F_1 \{\F_2 [(I_2^\b R_T  f)(a, \cdot)] (\eta) \})
(x'\eta)\nonumber\\
&=& (\F_1 \{|\eta|^{-\b}(\F f)(-a\eta, \eta)\})(x'\eta)\nonumber\\
&=&|\eta|^{-\b}\int_{\bbr^{m-1}} (\F f)(-a\eta, \eta)\, e^{ia \cdot
x'\eta }\, da\nonumber\\
&=&|\eta|^{1-\b-m}\int_{\bbr^{m-1}} (\F f)(w, \eta)\, e^{-iw \cdot
x'}\, dw\nonumber\\
&=& (2\pi)^{m-1} \,|\eta|^{1-\b-m}\, (\F_2[f(x',
\cdot)])(\eta).\nonumber\eea Hence, by (\ref{ash7}),
$$
(\F_2 I_2^\a [(\tilde R_T g)(x', \cdot)])(\eta)=(2\pi)^{m-1}
\,|\eta|^{1-\a-\b-m} (\F_2[f(x', \cdot)])(\eta),$$ which implies
 (\ref{cab}).
\end{proof}
Equality (\ref{cab}) gives a variety of pointwise inversion
formulas. Let $\bbd^\a_2=I_2^{-\a}$ be the corresponding partial
Riesz fractional derivative, which is well defined on functions $ f
\in \Phi(\bbr^m)$ in the Fourier terms: $(\F \bbd^\a_2
f)(y)=|y_m|^\a (\F f)(y)$.

\begin{corollary}\label {puu}  Let $\vp=R_T  f, \; f \in
\Phi(\bbr^m)$. The following pointwise inversion formulas are
contained in (\ref{cab}):

\be\label {ffo1}  \!\!\!\!\!\!\!\!\!\!f=(2\pi)^{1-m} \, \bbd_2^{m-1}
\tilde R_T \vp,\qquad \quad \qquad\quad(\b=0, \; \a =1-m);\quad\ee
\be\label {ffo2} \!\!\!\!\! \!\!\!\!\!f=(2\pi)^{1-m} \,  \tilde R_T
\bbd_2^{m-1}\vp,\qquad \quad \qquad\quad(\a=0, \; \b =1-m);\quad\ee
\be\label {ffo3}  f\!=\!(2\pi)^{1-m} \, \bbd_2^{(m-1)/2} \tilde R_T
\bbd_2^{(m-1)/2} \vp,\qquad (\a\!=\!\b\!=\!(1\!-\!m)/2);\ee
\be\label {ffo4}\!\!\!\!\!  \!\!\!\!\!f=(2\pi)^{1-m} \, \bbd_2^{m}
\tilde R_T I_2^1\vp,\qquad \quad \qquad\quad(\a =-m, \; \b
=1).\qquad\ee
\end{corollary}

The Riesz derivative $\bbd_2$, corresponding to $|y_m|$ in the
Fourier terms, can be expressed through the usual derivative
$\partial_m$ by the formula $\bbd_2=iH_2\partial_m$, where $H_2$
stands for the Hilbert transform in the last variable: \be\label
{ffo14} (H_2 \vp)(a,b)\!=\!\mbox {\rm p.v.}\frac{1}{\pi
i}\int_{\bbr}\frac{\vp(a,b_1)}{b-b_1}\, db_1, \quad (\F H_2
\vp)(y)\!=\!\sgn y_m \, (\F \vp)(y);\nonumber\ee see, e.g., \cite
{Graf}. In particular, for $m\!=\!2n\!+\!1, \; n\!\in \!\bbn$,  we
have the following
\begin{corollary}  \label {puu1} Let $\vp=R_T  f, \; f \in
\Phi(\bbr^{2n+1})$. Then
\bea
\label {ffo11}f&=&(2\pi)^{-2n} (-1)^n\,\partial_{2n+1}^{2n}\tilde R_T \vp\\
\label {ffo12}&=&(2\pi)^{-2n} (-1)^n\,\tilde R_T \partial_{2n+1}^{2n}\vp\\
\label {ffo13}&=&(2\pi)^{-2n}(-1)^n\,\partial_{2n+1}^{n}\tilde R_T \partial_{2n+1}^{n}\vp.\eea
\end{corollary}

In fact, (\ref{ffo12}) and (\ref{ffo12}) follow from (\ref{ffo11})
 because $\tilde R_T$ commutes with differentiation in the last variable.

\begin{remark} Inversion formula (\ref{roz1}) and those in  Corollaries
\ref{puu} and \ref{puu1} are pointwise analogues for $R_T$ of the
Strichartz's formulas (\ref{Siv}). In the next section we adapt
these formulas for the Radon transform on the Heisenberg group.
\end{remark}

\section{Connection between  the Radon
transforms $R, R_T, R_H$}

We denote $$ x\!=\!(x_1, \ldots,
x_m)\!=\!(x', x_m)\!\in\! \bbr^m, \qquad \theta\!=\!(\theta_1,
\ldots,\theta_m)\!=\!(\theta', \theta_m )\!\in \!S^{m-1},$$ and let
$h$ be a hyperplane which is not parallel to the $x_m$-axis. The
 following two parameterizations are available:
 \be
\label{par1} h=\{x: x_m=a\cdot x' +b\}, \qquad (a,b)\in \tilde
\bbr^m =\bbr^{m-1}\times \bbr,\ee \be \label{par2} h=\{x:
\theta\cdot x=t\},\qquad (\theta, t)\in \bbp^m=S^{m-1}\times \bbr,
\quad \theta_m \neq 0. \ee Different sets of parameters are related
by
 \be \label{change1}
 \theta=\frac{a-e_m}{\sqrt{1+|a|^2}}, \qquad
t=-\frac{b}{\sqrt{1+|a|^2}},\ee \be \label{change2}
a=-\frac{\theta'}{\theta_m}, \qquad b=\frac{t}{\theta_m}.\ee The
corresponding {\it transference operators} acting on functions $\vp
(a,b)$ and $\psi (\theta, t)$ have the form \be (\T\vp)(\theta,
t)=\vp \left (-\frac{\theta'}{\theta_m},\frac{t}{\theta_m}\right),
\ee \be\label{tinv} (\T^{-1} \psi)(a,b)=\psi \left
(\frac{a-e_m}{\sqrt{1+|a|^2}}, -\frac{b}{\sqrt{1+|a|^2}}\right).\ee

\begin{lemma} The following relations hold:
 \be \label{45}(Rf)(\theta,
t)=|\theta_m |^{-1} (\T R_T f)(\theta, t), \qquad  \theta_m \neq
0,\ee \be\label{46} (R_T f)(a,b)=(1+|a|^2)^{-1/2} (\T^{-1}
Rf)(a,b).\ee If $\psi (\theta, t)=\psi (-\theta, -t)$ for all
$(\theta, t)\in \bbp^m$, then for the dual transforms  (\ref{durt})
and (\ref{du1}) we have \be\label{47}  (R^*
\psi)(x)=2(\stackrel{*}{R_T} \T^{-1}\psi)(x).\ee
Furthermore,\be\label{48}  (\stackrel{*}{R_T}
\vp)(x)=\frac{1}{2}\,(R^*\T\vp)(x).\ee

\end{lemma}
\begin{proof} Equalities (\ref{45}) and (\ref{46}) follow from
(\ref{con1}). Equality (\ref{47}) can be obtained by making use of
(\ref {teq1}). Namely, \bea (R^*
\psi)(x)&=&\int_{S^{m-1}}\psi(\theta, \theta\cdot x)\,
d\theta\nonumber\\&=& \int_{\bbr^{m-1}}\Big[\psi \Big
(\frac{a+e_m}{\sqrt{1+|a|^2}}, \frac{(a+e_m)\cdot
x}{\sqrt{1+|a|^2}}\Big )\nonumber\\
&+&\!\psi \Big (\frac{a-e_m}{\sqrt{1+|a|^2}}, \frac{(a-e_m)\cdot
x}{\sqrt{1+|a|^2}}\Big )\Big ] \tilde d  a\nonumber\\&=&
2\int_{\bbr^{m-1}}\psi \left (\frac{a-e_m}{\sqrt{1+|a|^2}},
\frac{(a-e_m)\cdot x}{\sqrt{1+|a|^2}}\right )\tilde d
a\nonumber\\&=& 2\int_{\bbr^{m-1}} (\T^{-1}\psi)(a, x_m- a\cdot
x')\, \tilde d a=2(\stackrel{*}{R_T} \T^{-1}\psi)(x).\nonumber\eea
 Equality (\ref{48}) is a consequence of  (\ref{47}). It can also be obtained
directly by using (\ref{teq2}).
\end{proof}

Connection between  the Radon transforms $R_T$ and $R_H$ is given
in (\ref{ccon1a}).
 Let us present this equality  and a similar one for
 the dual transforms in the operator form. We define
 \bea
(\stackrel{*}{R_H}\vp)(\z,\t)&=&2\int_{\bbc^n}\vp((\z,\t) \circ (z,
0))\, \frac{dz}{(4+|\z+z|^2)^{n+1/2}}
\nonumber\\
&=&\label{dualRH}\int_{\bbc^n}
 \!\!\vp(z, \t\!-\!\frac{1}{2} \,Im\, (\z\cdot \bar z))\,
 \tilde{d}z, \quad \tilde{d}z=\frac{2dz}{(4\!+\!|z|^2)^{n+1/2}},\eea
and set
$$
x=(x_{(1)}, x_{(2)}, x_{2n+1}), \quad x_{(1)}=(x_1, \ldots, x_n), \quad x_{(2)}=(x_{n+1}, \ldots, x_{2n}),
$$
$$
a=(a_{(1)}, a_{(2)}), \quad a_{(1)}=(a_1, \ldots, a_n), \quad a_{(2)}=(a_{n+1}, \ldots, a_{2n}).
$$

\begin{lemma} \label{rath} Given functions $f$ and $\vp$ on $\bbh_n$, let \be(\Q
f)(x)\!=\!f(x_{(1)}\!+\!ix_{(2)}, x_{2n+1}), \quad (\tilde\Q
\vp)(a,b)\!=\!\vp (2a_{(2)}\!-\!2ia_{(1)}, b).\ee
Then
\be \label{rela1} \tilde\Q R_H=R_T \Q.\ee \be\label{rela2}
\Q\stackrel{*}R_H=\stackrel{*}{R_T}\tilde\Q.\ee The duality relation
\be\label{rela3}
\int_{\bbh_n}f(\z,\t)\,(\stackrel{*}{R_H}\vp)(\z,\t)\,d\z
d\t=\int_{\bbh_n}(R_Hf)(z,t)\,\vp(z,t)\,\tilde{d} z dt \ee holds provided
that at least one of these integrals is finite when $f$ and $\vp$
are replaced by $|f|$ and $|\vp|$.
\end{lemma}
\begin{proof} Equality (\ref{rela1}) follows from
(\ref{ccon1})-(\ref{ccon1a}). Furthermore, by (\ref{du1}),\bea
&&(\stackrel{*}{R_T}\tilde\Q \vp)(x)=\int_{\bbr^{2n}}(\tilde\Q
\vp)(a, x_{2n+1}-
a\cdot x')\, \tilde d  a\nonumber\\
&&=\int_{\bbr^{2n}}\vp (2a_{(2)}-2ia_{(1)}, x_{2n+1}- a\cdot
x')\,\frac{da}{(1+|a|^2)^{n+1/2}}
\nonumber\\
&&=2\int_{\bbr^{2n}}\vp (u+iv, x_{2n+1}- \frac{1}{2}(u\cdot x_{(2)}
-v\cdot x_{(1)}) \,\frac{du\, dv}{(4+|u|^2+|v|^2)^{n+1/2}}.
\nonumber\eea Setting $x_{2n+1}=\t, \;x_{(1)}=\xi, \; x_{(2)}=\eta,
\; \xi +i\eta=\z, \; u+iv=z$, we get \bea (\Q^{-1}\stackrel{*}R_H
\vp)(\z, \t)\!\!\!&=&\!\!\!2\!\int_{\bbr^{2n}}\!\!\!\vp (u\!+\!iv,
\t\!- \!\frac{1}{2}(u\cdot \eta \!-\!v\cdot \xi) \,\frac{du\,
dv}{(4\!+\!|u|^2\!+\!|v|^2)^{n+1/2}}\nonumber\\
&=&\int_{\bbc^{n}}\vp (z, \t- \frac{1}{2} \,Im\, (\z\cdot \bar z))\,
 \tilde{d}z.\nonumber\eea
This gives (\ref{rela2}).  The duality relation (\ref{rela3})
follows from (\ref{dual1}) by the same reasoning.
\end{proof}

\section{Inversion of the  Heisenberg Radon Transform}

Equalities (\ref{rela1}) and (\ref{rela2}) enable us to convert
inversion formulas for the transversal Radon transform $R_T$ (with
$m=2n+1$) into those for  the Heisenberg Radon transform $R_H$. We
skip routine calculations, which mimic the proof of Lemma
\ref{rath}.  As in (\ref{SRT}), we write $(\z, \t)$ for
 the argument of $f$ and $(z, t)$ for
 the argument of $R_H f$.

Let us start with Theorem \ref{441a} and set
 \be\label{ge1} g(\z, \t)=(2\pi)^{-2n}
\,(\stackrel{*}{R_H}\tilde\psi)(\z, \t),\ee $$ \tilde\psi
(z,t)=\sqrt{1+|z|^2/4}\, (R_H f)(z, t); $$
 \bea \label{diff1} \qquad (\Del^\ell_{(z,t)} g)(\z,\t)\!\!&=&\!\!\sum_{j=0}^\ell
{\ell \choose j} (-1)^j g(\z\!-\!jz, \t\!-\!jt), \quad \ell>2n,\\
\label{diff2} \qquad (\tilde \Del^k_{(z,t)}
g)(\z,\t)\!\!&=&\!\!\sum_{j=0}^k {k \choose j} (-1)^j
g(\z\!-\!\sqrt{j} z, \t\!-\!\sqrt{j} t),\quad k>n; \eea
$$
d_{n,\ell} \! = \! \int_{\bbr^{2n+1}}\!\! \frac{(1 \! - \!
e^{iy_1})^\ell}{|y|^{4n+1}} dy, \qquad   \tilde
d_{n,k}=\frac{2^{-2n}\pi^{n+1/2}}{\Gam (2n+1/2)}\intl_0^\infty
\frac{(1\!-\!e^{-t})^k}{t^{n+1}}\, dt.
$$
\begin{theorem} \label{inveH1} Let $ f \in L^p(\bbh^n), \; 1 \le p
<1+1/2n$. Then $f$ can be reconstructed by the formulas
\be\label{forH1} f(\z,\t) \! = \! \frac{1}{d_{n,\ell}}\int_{\bbh_n}
\frac{(\Del^\ell_{(z,t)} g)(\z,\t)}{(|z|^2+t^2)^{2n+1/2}} \,dz dt
 \nonumber\ee
 or
 \be\label{forH2} f(\z,\t) \! = \! \frac{1}{\tilde
d_{n,k}}\int_{\bbh_n} \frac{(\tilde\Del^k_{(z,t)}
g)(\z,\t)}{(|z|^2+t^2)^{2n+1/2}} \,dz dt, \nonumber\ee where
$\int_{\bbh_n}=\lim_{\varepsilon\rightarrow0}\int_{|z|^2+t^2>\varepsilon}$.
 The limit exists in the $L^p$-norm and in the a.e. sense. For
$f\in C_0(\bbh^n)\bigcap L^p(\bbh^n)$ it exists in the $\sup$-norm.
\end{theorem}

The next statement follows from Theorem \ref{inve2}.
\begin{theorem}  \label{inveH2} Let $f \in C^*_\del(\bbh_n)$,
$\del>2n$, and let $g$ be defined  by (\ref{ge1}).  Then
 \be \label{fH1}f(\xi +i\eta,\t)= (-\Del)^{n}\,g(\xi +i\eta,\t),\ee
 where
$$
\Del=\sum\limits_{k=1}^n\frac{\partial^2}{\partial
\xi_k^2}+\sum\limits_{k=1}^n\frac{\partial^2}{\partial
\eta_k^2}+\frac{\partial^2}{\partial \t^2}.$$
\end{theorem}

 Theorems \ref{441ab},  \ref{cbp},  and \ref{cbpx} can be
 reformulated in a similar way. We leave them to the interested
 reader.

Another series of pointwise inversion formulas for $R_H$ can be
derived in the framework of the corresponding Semyanistyi-Lizorkin
spaces.  Following Section \ref{sed}, we define \be  \label
{hfii}\Phi (\bbh_n)=\{\phi\in S(\bbh_n): \,\int_{\bbr} \phi(\z,
\t)\, \t^{k}
 \,d\t\!=\!0 \;\forall k\in \bbz_+, \; \forall \z \in\bbc^{n}\}.\ee

\begin{theorem} \label{hish}The Heisenberg  Radon transform $R_H$ is  an automorphism of the space $\Phi (\bbh_n)$. \end{theorem}
\begin{proof}
By (\ref{rela1}), $  R_H=\tilde\Q^{-1}R_T \Q$. Hence, Theorem \ref{hish} is a consequence of Theorem \ref{ish}.
\end{proof}
\begin{theorem} \label{inve1} A function $f \in \Phi (\bbh_n)$ can be reconstructed from its Heisenberg  Radon transform $\vp=R_H f$ by the formula
\be \label{lpl} f(\xi+i\eta, \t)= \frac{1}{(2\pi)^{2n+1}}\int_{\bbr^{2n+1}} w (y)\, e^{-i y\cdot (\xi,\eta, \t)}\, dy,\ee where
$$
w (y)=\int_{\bbr} \vp\Big (-\frac{2y_{(2)}-2iy_{(1)}}{y_{2n+1}}, b\Big )\, e^{iby_{2n+1}}\, db, \qquad y=(y_{(1)},y_{(2)},y_{2n+1}).
$$
 \end{theorem}
\begin{proof} Replace $f$ by $\Q f$ in (\ref{roz1}) and make use of (\ref{rela1}). This gives
$$
 f=\Q^{-1}\F^{-1}\Lam^{-1}\F_2 \tilde\Q \vp,$$  which coincides with (\ref{lpl}).
\end{proof}

Inversion formulas below do not contain the Fourier transform.

\begin{theorem} \label{inve2} A function $f \in \Phi (\bbh_n)$ can be reconstructed from  $\vp=R_H f$ by the following formulas:
\bea \qquad f (\z,\t)&=& (-1)^n\,(4\pi)^{-2n}\,\partial_\t^{2n} (R_H \vp)(\z,\t)\\
&=&
(-1)^n\,(4\pi)^{-2n}\,\partial_\t^{2n} (R_H \partial_t^{2n}\vp)(\z,\t)\\
&=&(-1)^n\,(4\pi)^{-2n}\,\partial_\t^{n} (R_H
\partial_t^{n}\vp)(\z,\t).\eea
\end{theorem}
\begin{proof} Since $R_H$ commutes with differentiation in the last variable, it suffices to check the first formula.
We replace $f$ by $\Q f$ in (\ref{ffo11}) and make use of
(\ref{rela1}). This gives
$$
 f=(-1)^n\,(2\pi)^{-2n} \,\Q^{-1}\partial_{2n+1}^{2n}\tilde R_T \tilde \Q\vp.$$
By (\ref{bpro}),
\bea
 (\tilde R_T \tilde \Q\vp)(x)&=&\int_{\bbr^{2n}} (\tilde \Q\vp)(a, -x'\cdot a+ x_{2n+1})\,da\nonumber\\
&=&\int_{\bbr^{2n}} \vp(2a_{(2)}-2ia_{(1)}, -x'\cdot a+ x_{2n+1})\,da\nonumber\\
&=&4^{-n}\int_{\bbr^{2n}} \!\!\vp(u+iv, \frac{1}{2}(v\cdot
x_{(1)}-u\cdot x_{(2)})+ x_{2n+1})\,dudv. \nonumber\eea Setting
$\z=\xi +i\eta$, $z=u+iv$, we have \bea (\Q^{-1}\tilde R_T \tilde
\Q\vp)(\z,\t)&=&4^{-n}\int_{\bbr^{2n}}
\vp(u+iv, \t-\frac{1}{2}(u\cdot\eta - v\cdot \xi))\,dudv\nonumber\\
&=&4^{-n}\int_{\bbc^n}
 \!\!\vp(z, \t\!-\!\frac{1}{2} \,Im\, (\z\cdot \bar z))\,
 dz=4^{-n}(R_H \vp)(\z,\t).\nonumber\eea This gives
 the result.
\end{proof}

Theorems \ref{inve1} and \ref{inve2} give precise  meaning to
Strichartz's formulas in (\ref{Siv}).

\section {Mixed-norm estimates for $R_T$ and $R_H$}

For $1\le q,r <\infty$, we define the following spaces with mixed
norm: \bea L^{q, r} (\bbp^m) &&\!\!\!\!\!\!\!= \{ \vp (\theta, t): \nonumber\\
&&\!\!\!\!\!\| \vp; \bbp^m\|_{q, r} \, = \, \left(\, \int_{S^{m-1}}
\left[\int_\bbr |\vp(\theta, t)|^r dt\right]^{q/r}
d\theta\right)^{1/q}<\infty\};\nonumber\eea

\bea L^{q, r} (\tilde\bbr^m) &&\!\!\!\!\!\!\!= \{ \psi (a,b): \nonumber\\
&&\!\!\!\!\! \| \psi; \tilde\bbr^m\|_{q, r} \, = \, \left(\,
\int_{\bbr^{m-1}} \left[\int_\bbr |\psi (a,b)|^r db\right]^{q/r}
da\right)^{1/q}<\infty\};\nonumber\eea

\bea \!\!L^{q, r} (\bbh_n) &&\!\!\!\!\!\!\!= \{ g (z, t):  \nonumber\\
&&\!\!\!\!\!\| g; \bbh_n\|_{q, r} \, = \, \left(\, \int_{\bbc^{n}}
\left[\int_\bbr |g (z, t)|^r dt\right]^{q/r}
dz\right)^{1/q}<\infty\}. \nonumber\eea

\begin{theorem} \label {OS} For $m\ge 2$ an
a priori
 inequality
$$\| Rf;  \bbp^m\|_{q, r} \le c_{p, q, r}  \| f \|_p$$
holds if and only if $1\le p < m/(m-1), \; q \le p^{\prime}\; \,(1/p
+ 1/p^{\prime} = 1)$ and
 $1/r = 1-m/p'$.
\end{theorem}

This well-known statement, which is due to  Oberlin and Stein
\cite{OS} and also  Strichartz \cite{Str1},
 implies the following result for the transversal Radon
transform.

\begin{theorem} \label {TRT} For $m\ge 2$ an
a priori
 inequality
$$\| R_Tf; \tilde\bbr^m\|_{q, r} \le c_{p, q, r}  \| f \|_p$$
holds if and only if \be\label{TRTa}1\le p < m/(m-1), \qquad q
=p^{\prime},\qquad 1/r = 1-m/p'.\ee
\end{theorem}
\begin{proof} Note that, unlike Theorem \ref{OS}, now we have a strict equality for $q$.
 Using (\ref{46}), (\ref{tinv}),   (\ref{teq2}), and
taking into account that $$(Rf)(-\theta,-t)=(Rf)(\theta,t),$$ we
obtain
\begin{eqnarray*}
&&\| R_Tf; \tilde\bbr^m\|_{q, r}^{q}\\
&&=\int_{\RR^{m-1}}\Big [\int_{\RR}\Big|(1+|a|^2)^{-1/2}
(Rf)\Big(\frac{a-e_m}{\sqrt{1+|a|^2}},
-\frac{b}{\sqrt{1+|a|^2}}\Big)\Big|^r db \Big]^{q/r}da \\
&&=\int_{\RR^{m-1}}\Big [\int_{\RR}\Big|
(Rf)\Big(\frac{a-e_m}{\sqrt{1+|a|^2}},
t\Big)\Big|^r dt \Big]^{q/r} (1+|a|^2)^{(1-r)q/2r}da \\
&&=\frac{1}{2}\int_{S^{m-1}}\Big
[\int_{\RR}\left|(Rf)(\theta,t)\right|^r\,dt\Big]^{q/r}\,|\theta_m|^{-m+(r-1)q/r}\,d\theta.
\end{eqnarray*}
Now we set $q =p^{\prime}$, $1/r = 1-m/p'$, and apply Theorem
\ref{OS}. This gives
$$
\| R_Tf; \tilde\bbr^m\|_{q, r}=2^{-1/q}\| Rf; \bbp^m\|_{q, r}\le
2^{-1/q}c_{p, q, r}  \| f \|_p,$$ as desired.

The necessity of (\ref{TRTa}) can be proved using the homogeneity
argument, as in \cite{Str2}.
\end{proof}

The mixed-norm estimate for the Heisenberg Radon transform now
 follows by Lemma \ref{rath}.
\begin{theorem} \label {mne2} For $n\ge 1$, an
a priori
 inequality
$$\| R_Hf; \bbh_n\|_{q, r} \le c_{p, q, r}  \| f \|_p$$
holds if and only if \be\label{TRTah}1\le p <1+1/2n, \qquad q
=p^{\prime},\qquad 1/r = 1/p-2n/p'.\ee
\end{theorem}
\begin{proof} We have
\bea &&\| R_Hf; \bbh_n\|_{q, r}^q=\int_{\bbr^n \times \bbr^n}\Big
[\int_\bbr |(R_H f)(u+iv,t)|^r dt \Big]^{q/r}\, dudv \nonumber\\
&& \mbox{\rm (set $u=2a^{(2)},
\quad v=-2a^{(1)}, \quad a=(a^{(1)}, a^{(2)}), \quad t=b$)}\nonumber\\
&&=4^n \int_{\bbr^{2n}}\Big [\int_\bbr |(R_T \Q f)(a,b)|^r db
\Big]^{q/r}\, da. \nonumber\eea Hence, by Theorem \ref {TRT},
$$
\| R_Hf; \bbh_n\|_{q, r} \le 4^{n/q} \,c_{p, q, r} ||Qf||_{L^p
(\bbr^{2n+1})}= 4^{n/q} \,c_{p, q, r} ||f||_{L^p (\bbh_n)}.$$ This
gives the result.
\end{proof}

Theorems \ref{TRT} and \ref{mne2}  generalize Corollary 4.2 in
\cite{Str2} and extend it to the full range of $p$; cf. Theorem
\ref{t2} and discussion after it.


\begin{thebibliography}{[ASMR]}

\bibitem [ABC] {ABC} M.L. Agranovsky,  C.A. Berenstein, and D.C. Chang,
\textit{ Morera theorem for holomorphic $H\sp p$ spaces in the
Heisenberg group},  J. Reine Angew. Math.,  \textbf{443}  (1993),
49--89.

\bibitem [ABR] {ABR}  S. Axler,  P. Bourdon, and W. Ramey, Harmonic Function Theory, Springer, 1992.

\bibitem [Ber] {Ber} C. A. Berenstein, \textit{Radon transforms, wavelets, and
applications},  Integral geometry, Radon transforms and complex
analysis (Venice, 1996),  1--33, Lecture Notes in Math., 1684,
Springer, Berlin, 1998.

\bibitem [BCT] {BCT}   C. A. Berenstein, D.C. Chang, and J. Tie, Laguerre calculus and its applications on the Heisenberg group.
AMS/IP Studies in Advanced Mathematics, \textbf{22}. American
Mathematical Society, Providence, RI; International Press,
Somerville, MA, 2001.


\bibitem [C1] {C1} E. J. Cand\`es, Ridgelets: theory and
applications. Ph.D. Thesis, Technical Report, Department of
Statistics, Stanford University, 1998.

\bibitem [C2] {C2} \bysame, \textit{Harmonic analysis of neural networks},
Appl. Comput. Harmon. Anal. \textbf{6} (1999), no. 2, 197--218.

\bibitem [CK] {CK} M. Cowling,  and A. Korányi,  \textit{Harmonic analysis on Heisenberg type
groups from a geometric viewpoint},  Lie group representations, III
(College Park, Md., 1982/1983),  60--100, Lecture Notes in Math.,
1077, Springer, Berlin, 1984.


\bibitem [Dea] {Dea}  S. R. Deans, The Radon transform and some of its applications,
Dover Publ. Inc., Mineola, New York, 2007.


\bibitem  [Ehr] {Ehr} L. Ehrenpreis, The universality of the Radon transform, Oxford University Press,
2003.


\bibitem [Fe1]  {Fe1} R. Felix, \textit{Radon-Transformation auf nilpotenten Lie-Gruppen} (German)
[Radon transform on nilpotent Lie groups],  Invent. Math.  112
(1993),  no. 2, 413--443.


\bibitem [Fe2]  {Fe2} \bysame, \textit{A general approach to the Radon transform},
 Proceedings of the conference "Program Systems: Theory and Applications" in Pereslavl-Zalessky, Russia, October 2006
ISBN 3-8322-2147-6, 31 pp.

\bibitem [FJW] {FJW} M. Frazier, B. Jawerth,  and G. Weiss,
 Littlewood-Paley theory and the study of function spaces,
CBMS Reg. Conf. Ser. in Math.,  no.  79, Amer. Math. Soc., Providence,
R.I.,  1991.

\bibitem  [GGG] {GGG} I. M. Gelfand, S. G. Gindikin,  and M. I. Graev,  Selected
topics in integral geometry,  Translations of Mathematical
Monographs,
 220,    American Mathematical Society, Providence, RI, 2003.

\bibitem [GGV] {GGV} I. M. Gelfand,  M. I. Graev, and N. Ja.
Vilenkin,  Generalized Functions, Vol 5, Integral geometry and
representation theory, Academic Press,  1966.


 \bibitem  [GS1] {GS1} D. Geller,  and  E. M. Stein,  \textit{Singular convolution operators on the Heisenberg
 group},  Bull. Amer. Math. Soc. (N.S.)  \textbf{6}  (1982), no. 1, 99--103.

 \bibitem  [GS2] {GS2} \bysame, \textit{Estimates for singular convolution operators on the Heisenberg group},
   Math. Ann.  \textbf{267}  (1984),  no. 1, 1--15.

\bibitem  [GRy]{GRy} I. S. Gradshteyn,  and  I. M. Ryzhik, Table of integrals, series and products,
Academic Press, 1980.


\bibitem  [Graf] {Graf} L. Grafakos,  Classical Fourier Analysis,
Second Edition, Graduate Texts in Math., no \textbf{249}, Springer,
New York, 2008.

\bibitem [He1] {He1} J. He,  \textit{An inversion formula of the Radon transform on the
Heisenberg group}  Canad. Math. Bull. \textbf{47}  (2004), 389--397.

\bibitem [He2] {He2} \bysame, \textit{A characterization of inverse Radon transform on the Laguerre
hypergroup},  J. Math. Anal. Appl. \textbf{318}  (2006), 387--395.


\bibitem [HL1] {HL1} J. He,  and H. Liu,   \textit{Inversion of the Radon
transform associated with the classical domain of type one},
Internat. J. Math.  16  (2005),   875--887.


\bibitem [HL2] {HL2}  \bysame, \textit{Admissible wavelets and
inverse Radon transform associated with the affine homogeneous
Siegel domains of type II}, Comm. Anal. Geom., \textbf {15} (2007),
1-28.

\bibitem [Hel] {Hel} S. Helgason,   The Radon transform,
Birkh\"auser, Boston, Second edition, 1999.

\bibitem [KR] {KR} A. Kaplan,  and F. Ricci,  \textit{Harmonic analysis on groups of Heisenberg
type}.  Harmonic analysis (Cortona, 1982),  416--435, Lecture Notes
in Math., \textbf{992}, Springer, Berlin, 1983.

\bibitem [Kat] {Kat}   A. I. Katsevich,  \textit{Range of the Radon transform on functions which
do not decay fast at infinity}, SIAM J. Math. Anal.,
 \textbf{28} (1997),  852--866.

\bibitem [Kor] {Kor}   A. Kor\'anyi,  \textit{Geometric aspects of analysis on the Heisenberg group},  Topics
    in modern harmonic analysis, Vol. I, II (Turin/Milan, 1982),  209--258, Ist. Naz. Alta Mat. Francesco Severi, Rome, 1983.

\bibitem [LH] {LH} P. Liu, and J. He,  \textit{Inversion of the Radon transform on the
product Laguerre hypergroup by using generalized wavelets},  Int. J.
Comput. Math.  \textbf{84}  (2007),  no. 3, 287--295.

\bibitem [Liz1] {Liz1} P.I. Lizorkin,  \textit{Generalized Liouville differentiation and
functional spaces $L_p^r(E_n)$. Imbedding theorems}, Matem. Sb.,
 {\bf 60}(120) (1963), 325--353 (Russian).


\bibitem [Liz2] {Liz2}\bysame,  \textit{Generalized Liouville differentiation and the
method of multipliers in the theory of imbeddings of classes of
differentiable functions},  Proc. Steklov Inst. Math., {\bf 105}
(1969), 105--202.


\bibitem [Liz3] {Liz3} \bysame,  \textit{Operators connected with fractional differentiation and classes of
differentiable functions},  Proc. Steklov Inst. Math., {\bf 117}
(1972), 251--286.


\bibitem [Mar] {Mar} A. Markoe,  Analytic tomography,  Encyclopedia of Mathematics and its
Applications {\bf  106},  Cambridge Univ. Press,  2006.

\bibitem [Mi] {Mi} S. G. Mikhlin,  Mathematical physics, an advanced
course, North-Holland Publ. Company, Amsterdam, 1970.

\bibitem [MRS] {MRS} D. M\"uller, F. Ricci, and E.M. Stein,\textit{Marcinkiewicz
multipliers and multi-parameter structure on Heisenberg (-type)
groups. I.},  Invent. Math.,  {\bf 119}  (1995),   199--233.

\bibitem [Na] {Na} F. Natterer, The mathematics of computerized tomography. Wiley,
New York, 1986.

\bibitem [NT] {NT}  M. M. Nessibi,  and K. Trim\`eche,   \textit{Inversion of the Radon transform on the
Laguerre hypergroup by using generalized wavelets}, J. Math. Anal.
Appl., \textbf{208} (1997), 337--363.

\bibitem [NS] {NS} V. A. Nogin,  and S. G. Samko,  \textit{Some applications of potentials and
approximative inverse operators in multi-dimensional fractional
calculus}, Fractional Calculus \& Applied Analysis \textbf{2}
(1999), no 2, 205-228.

\bibitem [OS] {OS} D. M. Oberlin,  and E. M. Stein,   \textit{Mapping properties of the
Radon transform},  Indiana Univ. Math. J., {\bf 31} (1982),
641--650.


\bibitem [OR] {OR}  E. Ournycheva,  and  B. Rubin, \textit{Semyanistyi's integrals and Radon
transforms on matrix spaces}, The Journal
 of Fourier Analysis and Applications, {\bf 14} (2008),  60--88.

\bibitem  [Pal] {Pal}  V. Palamodov, Reconstructive integral geometry, Monographs in
Mathematics, 98. Birkhäuser Verlag, Basel, 2004.

\bibitem [PZ] {PZ} L. Peng, and G. Zhang,\textit{ Radon transform on H-type and
Siegel-type nilpotent groups},  Internat. J. Math.,   \textbf{ 18 }
(2007), 1061--1070.


\bibitem [QCK] {QCK} E.T. Quinto,  M.  Cheney, P. Kuchment (Editors), Tomography, impedance
imaging, and integral geometry, Lectures in Applied Mathematics, 30,
American Mathematical Society, Providence, RI, 1994.

\bibitem [Rad]{Rad}  J. Radon, \textit{\"Uber die Bestimmung von
Funktionen durch ihre Integralwerte l\"angs gewisser
Mannigfaltigkeiten}, Ber. Verh. S\"achs. Akad. Wiss. Leipzig, Math.
- Nat. Kl., \textbf{69} (1917), 262--277 (English translation in
\cite{Dea}).

\bibitem [RK] {RK}  A. G. Ramm,  and  A. I. Katsevich,    The Radon transform and local
tomography, CRC Press, Boca Raton, 1996.

\bibitem [Rou]{Rou}  F. Rouvi\`ere, \textit{Inverting Radon transforms:
The group-theoretic approach}, L'Enseignement Math.  \textbf{47}
(2001), 205-252.

\bibitem [Ru1] {Ru1} B. Rubin, Fractional integrals and
potentials,
  Pitman Monographs and Surveys in Pure and Applied Mathematics,
  \textbf{82},  Longman, Harlow, 1996.

  \bibitem [Ru2]{Ru2} \bysame,  \textit{Fractional calculus and wavelet transforms in
integral geometry}, Fractional Calculus and Applied Analysis
\textbf{1} (1998), 193--219.

\bibitem [Ru3] {Ru3} \bysame, \textit{Spherical Radon transform
and related wavelet transforms}, Appl. and Comp. Harmonic Anal.
\textbf{5} (1998), 202-215.

\bibitem [Ru4]{Ru4} \bysame, \textit{ Inversion formulas for the spherical
 Radon transform and the generalized cosine transform},  Advances in Appl.
Math. \textbf{29} (2002), 471--497.

\bibitem [Ru5]{Ru5} \bysame,  \textit{ Radon, cosine, and sine transforms on real hyperbolic space},
 Advances in Math. \textbf{170} (2002), 206--223.

\bibitem [Ru6] {Ru6} \bysame, \textit{Notes on Radon transforms in integral geometry}, Fractional
Calculus and Applied Analysis  \textbf{6} (2003), 25-72.

\bibitem [Ru7] {Ru7} \bysame, \textit{ Reconstruction of functions from their integrals over $k$-planes},
Israel J. of Math., {\bf 141} (2004), 93-117.


\bibitem [Ru8] {Ru8} \bysame, \textit{Convolution-backprojection method for the $k$-plane transform,
and Calder\'on's identity for ridgelet transforms}, Appl. Comput.
Harmon. Anal. {\bf 16} (2004), 231-242.


\bibitem [RZ]{RZ}  B. Rubin,  and G. Zhang,   \textit{Generalizations of the
Busemann-Petty problem for sections of convex bodies}, J. Funct.
Anal. \textbf{213} (2004),  473--501.


 \bibitem [Sam1] {Sam1} S. G. Samko,  \textit{Test functions vanishing on a given set, and division by a function},
  Mat. Zametki, {\bf 21} (1977), No. 5, 677--689 (Russian).

\bibitem [Sam2] {Sam2}  \bysame, Denseness of Lizorkin-type spaces $\Phi_V$ in $L_p(R^n)$,
{\it Mat. Zametki}, {\bf 31} (1982), No. 6, 655--665 (Russian).

\bibitem [Sam3] {Sam3}   \bysame, Hypersingular integrals and their applications,
 Taylor \& Francis, Series: Analytical Methods and Special Functions, Vol. 5, 2002.

\bibitem [SKM] {SKM}   S. G. Samko,  A. A. Kilbas, and O. I.
Marichev,    Fractional integrals and derivatives. Theory and
applications, Gordon and Breach Sc. Publ., New York, 1993.

\bibitem [SS] {SS} N. G. Samko,  and S. G. Samko,  \textit{
On approximate definition of  fractional differentiation},
 Fractional Calculus and Applied Analysis \textbf{2} (1999),  329--342.


\bibitem [Se1] {Se1}  V.I. Semyanistyi, \textit{On some integral transformations
in Euclidean space}, Dokl. Akad. Nauk SSSR {\bf 134} (1960), 536-539
(Russian).

\bibitem [Se2] {Se2}  \bysame, \textit{ Homogeneous functions and some problems of integral geomery
in spaces of constant cuvature}, Sov. Math. Dokl.   \textbf{2}
(1961),  59-61.


\bibitem [SSW] {SSW} K. T. Smith, D. C. Solmon,  and S. L. Wagner,
\textit{Practical and mathematical aspects of the problem of
reconstructing objects from radiographs}, Bull. of Amer. Math. Soc.
\textbf{83} (1997), 1227-1270.

\bibitem [So] {So} D. C. Solmon,   \textit{A note on $k$-plane integral
transforms},  J. Math. Anal. Appl. \textbf {71} (1979), 351--358.


\bibitem [St1] {St1} E. M. Stein,   Singular integrals and
differentiability properties of functions. Princeton Univ. Press,
Princeton, NJ. 1970.

\bibitem [St2] {St2} \bysame, Harmonic analysis, real variable methods, orthogonality, and
oscillation integrals, Princeton Univ. Press, Princeton, NJ, 1993.



\bibitem [Str1] {Str1} R. S. Strichartz,    \textit{ $L^p$-estimates for Radon transforms in
Euclidean  and non-euclidean spaces},  Duke Math. J. \textbf{48}
(1981), 699--727.


\bibitem [Str2] {Str2}  \bysame,  \textit{$L\sp p$ harmonic analysis and Radon transforms on the Heisenberg
group},  J. Funct. Anal.  {\bf 96} (1991),  350--406.


\bibitem [Than] {Than} S. Thangavelu,  Harmonic analysis on the Heisenberg group.
Progress in Mathematics, 159. Birkhäuser Boston, Inc., Boston, MA,
1998.


\bibitem [Vl] {Vl} V. S. Vladimirov,  The equations of mathematical physics,
 ``Nauka'', Moscow, 1988  (in Russian).


\end{thebibliography}
\end{document}